# ASYMPTOTIC NORMALITY OF THE $K$-CORE IN RANDOM GRAPHS

BY SVANTE JANSON AND MALWINA J. LUCZAK[1]

*Uppsala University and London School of Economics*

We study the $k$-core of a random (multi)graph on $n$ vertices with a given degree sequence. In our previous paper [*Random Structures Algorithms* **30** (2007) 50–62] we used properties of empirical distributions of independent random variables to give a simple proof of the fact that the size of the giant $k$-core obeys a law of large numbers as $n \to \infty$. Here we develop the method further and show that the fluctuations around the deterministic limit converge to a Gaussian law above and near the threshold, and to a non-normal law at the threshold. Further, we determine precisely the location of the phase transition window for the emergence of a giant $k$-core. Hence, we deduce corresponding results for the $k$-core in $G(n,p)$ and $G(n,m)$.

**1. Introduction.** Let $k \geq 2$ be an integer, fixed throughout the paper. The $k$-core of a graph $G$ is the largest induced subgraph of $G$ with minimum vertex degree at least $k$. (It is easy to see that this is well defined, but note that the $k$-core may be empty. It is customary to say that a $k$-core exists if it is nonempty.)

The question whether or not a nonempty $k$-core exists in a random graph, together with questions concerning the size and structure of the $k$-core when it does exist, have attracted a great deal of attention over the recent years. Starting with the pioneering papers by Bollobás [3] and Łuczak [21], many authors have studied various types of random graphs and also hypergraphs. A milestone was the paper of Pittel, Spencer and Wormald [27], who found the threshold for the appearance of a nonempty $k$-core in the random graphs $G(n,p)$ and $G(n,m)$, as well as the size of the nonempty $k$-core. Several different proofs of this result have been given since—see our own proof in [17] and the references therein. There have also been a number of papers

Received December 2006; revised December 2006.
[1]Supported in part by the Nuffield Foundation.

studying the $k$-core of a random graph with a specified degree sequence, for example, Fernholz and Ramachandran [8, 9], Cooper [6], Molloy [23] (the last two references also consider hypergraphs) and Janson and Luczak [17]. Further, related models of random graphs have been studied by Kim [20], Darling and Norris [7] and Cain and Wormald [5]. [Many of the papers listed above deduce results for $G(n,p)$ or $G(n,m)$ from their main result.]

In our previous paper [17] we showed that a version of the standard core-finding algorithm leads to simple proofs of the results on the existence and size of the $k$-core in random graphs with a prescribed degree sequence (under certain conditions), and hence, also in the random graphs $G(n,p)$ and $G(n,m)$. Our main probabilistic tools are standard results on the convergence of empirical distributions of independent random variables, applied to balls-and-bins processes associated with the algorithm. In the present paper we show that a more refined study of these processes leads to Gaussian limit laws for the processes, and that these in turn imply a Gaussian limit law for the size of the $k$-core away from the threshold, as well as a precise description of the threshold and the size of the threshold $k$-core. In particular, we show that, for $G(n,m)$, the width of the threshold is of the order $\sqrt{n}$ edges.

Given a graph $G$, let $v(G)$ and $e(G)$ denote the sizes of its vertex and edge sets respectively. We assume that $v(G) = n$ and consider asymptotics as $n \to \infty$. We say that an event holds whp (*with high probability*), if it holds with probability tending to 1 as $n \to \infty$. All unspecified limits in this paper are as $n \to \infty$. We use $O_{\mathrm{p}}$ and $o_{\mathrm{p}}$ in the standard way (see, e.g., Janson, Łuczak and Ruciński [18]); for example, if $(X_n)$ is a sequence of random variables, then $X_n = O_{\mathrm{p}}(1)$ means "$X_n$ is bounded in probability" and $X_n = o_{\mathrm{p}}(1)$ means that $X_n \xrightarrow{\mathrm{P}} 0$.

Let us first introduce some notation and recall the main result of Pittel, Spencer and Wormald [27].

For $\mu \geq 0$, let $\mathrm{Po}(\mu)$ denote a Poisson random variable with mean $\mu$. We denote the Poisson probabilities by

$$\pi_j(\mu) := \mathbb{P}(\mathrm{Po}(\mu) = j) = \mu^j e^{-\mu}/j!, \tag{1.1}$$

for integer $j \geq 0$, and let

$$\psi_j(\mu) := \mathbb{P}(\mathrm{Po}(\mu) \geq j) = \sum_{i=j}^{\infty} \pi_i(\mu). \tag{1.2}$$

Note that, for $j \geq 1$, $\pi'_j(\mu) = \pi_{j-1}(\mu) - \pi_j(\mu)$ and $\psi'_j(\mu) = \pi_{j-1}(\mu)$. Also, let

$$c_k := \inf_{\mu > 0} \mu/\psi_{k-1}(\mu).$$

For $\lambda > c_k$, or $\lambda = c_k$ and $k \geq 3$, we use $\mu_k(\lambda) > 0$ to denote the largest solution to $\mu/\psi_{k-1}(\mu) = \lambda$.

Pittel, Spencer and Wormald [27] discovered that $p = c_k/n$ is the threshold for the appearance of a nonempty $k$-core in the graph $G(n, p)$; equivalently, $m = c_k n/2$ is the threshold in the graph $G(n, m)$. More precisely, their main result is the following. [We write $G(n, \lambda_n/n)$ for the random graph $G(n, p)$ with $p = \lambda_n/n$; we continue to write $G(n, m)$, but we will only consider the corresponding case $m = m(n) = \lambda_n n/2$, and we will assume $\lambda_n \to \lambda < \infty$.]

THEOREM 1.1 (Pittel, Spencer and Wormald [27]).  *Consider the random graph $G(n, \lambda_n/n)$, where $\lambda_n \to \lambda \geq 0$. Let $k \geq 2$ be fixed and let* $\text{Core}_k = \text{Core}_k(n, \lambda_n)$ *be the $k$-core of $G(n, \lambda_n/n)$:*

(i) *If $\lambda < c_k$ and $k \geq 3$, then* $\text{Core}_k$ *is empty whp.*
(ii) *If $\lambda > c_k$, then* $\text{Core}_k$ *is nonempty whp, and $v(\text{Core}_k)/n \xrightarrow{\text{p}} \psi_k(\mu_k(\lambda))$, $e(\text{Core}_k)/n \xrightarrow{\text{p}} \mu_k(\lambda)\psi_{k-1}(\mu_k(\lambda))/2 = \mu_k(\lambda)^2/(2\lambda)$.*

*The same results hold for the random graph $G(n, m)$, for any sequence $m = m(n)$ with $\lambda_n := 2m/n \to \lambda$.*

Note that part (i) does not hold for $k = 2$. In fact, $c_2 = 1$, but for any $\lambda > 0$ there is a positive limiting probability that there are cycles and thus a nonempty 2-core. Nevertheless, if $\lambda < c_2$, the core is small, $v(\text{Core}_k) \leq e(\text{Core}_k) = O_p(1)$ and, in particular, $v(\text{Core}_k)/n \xrightarrow{\text{p}} 0$.

Pittel, Spencer and Wormald [27] further obtained some refinements of this result. In particular, they found that in the case $\lambda = c_k$, for any $\delta < 1/2$, (i) applies when $\lambda_n < c_k - n^{-\delta}$ and (ii) applies when $\lambda_n > c_k + n^{-\delta}$.

We now state our main theorems. The first shows that above the threshold the numbers of vertices and edges in the $k$-core have asymptotic normal distributions.

THEOREM 1.2.  *Consider either $G(n, \lambda_n/n)$ or $G(n, m)$, where $m = \lambda_n n/2$, and let $\text{Core}_k$ be the $k$-core of this random graph. Assume that $k \geq 2$ and that $\lambda_n \to \lambda > c_k$. Then, with $\hat{\mu}_n := \mu_k(\lambda_n)$,*

$$n^{-1/2}(v(\text{Core}_k) - \psi_k(\hat{\mu}_n)n, e(\text{Core}_k) - \tfrac{1}{2}\hat{\mu}_n\psi_{k-1}(\hat{\mu}_n)n) \xrightarrow{\text{d}} (Z_v, Z_e),$$

*where $Z_v$ and $Z_e$ are jointly Gaussian random variables with mean 0 and a nonsingular covariance matrix. More precisely, let $\hat{\mu} := \mu_k(\lambda)$ [so that $\hat{\mu} = \lambda\psi_{k-1}(\hat{\mu})$], and define*

(1.3)     $a_v := \pi_{k-1}(\hat{\mu})/(\psi_{k-1}(\hat{\mu}) - \hat{\mu}\pi_{k-2}(\hat{\mu}))$,

(1.4)     $a_e := (\psi_{k-1}(\hat{\mu}) + \hat{\mu}\pi_{k-2}(\hat{\mu}))/(\psi_{k-1}(\hat{\mu}) - \hat{\mu}\pi_{k-2}(\hat{\mu}))$,

(1.5)     $\widehat{\sigma}_{\nu\varkappa} := \sigma_{\nu\varkappa}(\hat{\mu}/\lambda) + \sigma^*_{\nu\varkappa}(\hat{\mu}/\lambda), \qquad \nu, \varkappa \in \{B, H, L\}$,

where $\sigma_{\nu\varkappa}(p) = \sigma_{\nu\varkappa}(p; \mathrm{Po}(\lambda))$ are as in Theorem 3.1 given by (5.37)–(5.47) and $\sigma_{\nu\varkappa}^*(p) = \sigma_{\nu\varkappa}^*(p; \lambda)$ are as in Theorem 8.5 [different for $G(n,p)$ and $G(n,m)$], then

$$(1.6) \qquad \mathrm{Var}(Z_v) = \widehat{\sigma}_{BB} - 2a_v\widehat{\sigma}_{BL} + a_v^2\widehat{\sigma}_{LL},$$

$$(1.7) \qquad \mathrm{Var}(Z_e) = \tfrac{1}{4}(\widehat{\sigma}_{HH} - 2a_e\widehat{\sigma}_{HL} + a_e^2\widehat{\sigma}_{LL}),$$

$$(1.8) \qquad \mathrm{Cov}(Z_v, Z_e) = \tfrac{1}{2}(\widehat{\sigma}_{BH} - a_e\widehat{\sigma}_{BL} - a_v\widehat{\sigma}_{HL} + a_v a_e\widehat{\sigma}_{LL}).$$

At the threshold, we have the following companion result.

THEOREM 1.3. *Consider either $G(n, \lambda_n/n)$ or $G(n,m)$, where $m = \lambda_n n/2$, and let $\mathrm{Core}_k$ be the $k$-core of this random graph. Assume that $k \geq 3$ and that $\lambda_n \to c_k$. Let $\hat{\mu} := \mu_k(c_k)$, $\hat{p} := \hat{\mu}/c_k = \psi_{k-1}(\hat{\mu})$, and*

$$(1.9) \qquad \widehat{\beta} := (\hat{\mu} - k + 2)\pi_{k-2}(\hat{\mu}) > 0.$$

*Let further $\sigma^2 := \sigma_{LL}(\hat{p}) + \sigma_{LL}^*(\hat{p}) > 0$, where $\sigma_{LL}(\hat{p}) = \sigma_{LL}(\hat{p}; \mathrm{Po}(c_k))$ is given by (5.47) and $\sigma_{LL}^*(p) = \sigma_{LL}^*(p; c_k)$ is as in Theorem 8.5. Then:*

  (i) *If $n^{1/2}(\lambda_n - c_k) \to -\infty$, then whp $\mathrm{Core}_k$ is empty.*
  (ii) *If $n^{1/2}(\lambda_n - c_k) \to \gamma \in (-\infty, \infty)$, then*

$$\mathbb{P}(\mathrm{Core}_k \neq \varnothing) \to \Phi(\hat{p}^2\gamma/\sigma),$$

*where $\Phi$ is the standard normal distribution function, and, with $Z \sim N(0,1)$,*

$$(n^{-3/4}(v(\mathrm{Core}_k) - \psi_k(\hat{\mu})n, e(\mathrm{Core}_k) - \tfrac{1}{2}\hat{\mu}\psi_{k-1}(\hat{\mu})n) \mid \mathrm{Core}_k \neq \varnothing)$$
$$\to ((2/\widehat{\beta})^{1/2}\sqrt{\sigma Z + \hat{p}^2\gamma}(\pi_{k-1}(\hat{\mu}), \hat{p}) \mid Z > -\hat{p}^2\gamma/\sigma).$$

  (iii) *If $n^{1/2}(\lambda_n - c_k) \to +\infty$, then whp $\mathrm{Core}_k$ is nonempty. Moreover, with $\hat{\mu}_n := \mu_k(\lambda_n)$,*

$$(\lambda_n - c_k)^{1/2}n^{-1/2}(v(\mathrm{Core}_k) - \psi_k(\hat{\mu}_n)n, e(\mathrm{Core}_k) - \tfrac{1}{2}\hat{\mu}_n\psi_{k-1}(\hat{\mu}_n)n)$$
$$\xrightarrow{\mathrm{d}} (\pi_{k-1}(\hat{\mu})Z', \hat{p}Z'),$$

*where $Z' \sim N(0, \sigma^2/(2\widehat{\beta}\hat{p}^2))$.*

In particular, Theorem 1.3 shows that the the width of the threshold for existence of a nonempty $k$-core is of the order $\sqrt{n}$ edges for $G(n,m)$ [and thus of the order $n^{-1/2}$ for $G(n,p)$, as a function of $p$], which improves the result $O(n^{1/2+\varepsilon})$ for every $\varepsilon > 0$ given by Pittel, Spencer and Wormald [27]. More precisely, we have the following simple corollary.

THEOREM 1.4. *Let $k \geq 3$. Start with the empty graph with $n$ isolated vertices and add edges at random, one by one, uniformly over all possible positions. Let $M$ be the number of edges when the graph first has a nonempty $k$-core. Then,*

$$n^{-1/2}\left(M - \frac{c_k}{2}n\right) \overset{\mathrm{d}}{\longrightarrow} N(0, \sigma_k^2),$$

*where $\sigma_k^2 := \sigma^2/(4\hat{p}^4)$ with $\sigma^2$ and $\hat{p}$ as in Theorem 1.3 for $G(n, m)$.*

Our formula for $\sigma_k^2$ is rather complicated. Numerical evaluations yield $\sigma_3^2 \approx 0.763$ and $\sigma_4^2 \approx 0.885$. (We have $c_3 \approx 3.35$ and $c_4 \approx 5.15$ as found by Pittel, Spencer and Wormald [27].)

REMARK 1.5. It follows from the above that the order of the typical random fluctuations of the number of vertices or edges in the $k$-core is $n^{1/2}$ above the threshold (Theorem 1.2), decreases like $n^{1/2}|\lambda_n - c_k|^{-1/2}$ as $\lambda_n$ approaches the threshold $c_k$ [Theorem 1.3(iii)], and becomes $n^{3/4}$ right at the threshold [Theorem 1.3(ii)]. Asymptotic normality holds above and near the threshold, but not at the threshold itself. The numbers of vertices and edges in the $k$-core (when nonempty) are asymptotically linearly dependent close to the threshold (Theorem 3.5), but not away from the threshold (Theorem 3.4).

The main idea in the proofs is to use the version of the core-finding process employed in [17] together with a martingale limit theorem by Jacod and Shiryaev [11]; we are then able to show that certain stochastic processes describing the progress of the core-finding algorithm are asymptotically Gaussian. Once this is done, we read off size of the $k$-core from these processes; this is rather straightforward although the details are time-consuming.

REMARK 1.6. Our method works for a fixed random graph, and is not directly applicable to studying the evolution of the $k$-core in the random graph process obtained by adding edges one by one at random. We can, however, say a little more. Consider $G(n, m)$ for two values of $m$, of the type $c_k n/2 + \gamma'_n n^{1/2}/2$ and $c_k n/2 + \gamma''_n n^{1/2}/2$ for two convergent sequences $\gamma'_n$ and $\gamma''_n$ with $\gamma'_n < \gamma''_n$. We suppose that the random graphs are coupled so that the second random graph is obtained by adding edges at random to the first. Then Theorem 1.3(ii) holds for each of the two random graphs separately, but it follows from the proof below, by coupling the processes of balls and bins used there in a straightforward way, that indeed Theorem 1.3(ii) holds jointly for both random graphs *with the same $Z$*.

The same applies to any finite number of stages in the evolution. In particular, we may consider $m_j(n) = \lfloor (c_k n + \gamma_j n^{1/2})/2 \rfloor$ with $\gamma_j = j/N$ for $j = -N^2, \ldots, N^2$, for a fixed $N$ at first (and letting $n \to \infty$). With probability $1 - \epsilon(N) + o(1)$, where $\varepsilon(N) = \mathbb{P}(|Z| > N\hat{p}^2/\sigma)$, the $k$-core will first appear somewhere inside this grid, say when $m_{j_0-1} < m \leq m_{j_0}$. More precisely, we may as in the proof below assume that the limit in Theorem 1.3(ii) holds a.s., and then this holds a.s. if $j_0$ is the smallest integer such that $\sigma Z + \hat{p}^2 \gamma_{j_0} > 0$. In this case, $\sigma Z + \hat{p}^2 \gamma_{j_0} \leq \hat{p}^2/N$, and it follows that the number of vertices in the $k$-core of $G(n, m_{j_0})$ is $\psi_k(\hat{\mu})n + O(n^{3/4}N^{-1/2})$. Hence, with probability $1 - \varepsilon(N) + o(1)$, the number of vertices in the first nonempty $k$-core is $\psi_k(\hat{\mu})n + O(n^{3/4}N^{-1/2})$.

Letting $N \to \infty$, we find that $\epsilon(N) \to 0$ (so that the probability that the $k$-core first appears within the grid approaches 1), and that the number of vertices in the first nonempty $k$-core is $\psi_k(\hat{\mu})n + o_p(n^{3/4})$. Hence, the random fluctuation is smaller than what we see in Theorem 1.3(ii) for a fixed $m$. (We cannot say how much smaller.) In other words, the order $n^{3/4}$ fluctuations must come from the fluctuations in the time the $k$-core appears. We also see that the $k$-core grows rapidly in the beginning due to the term $\sqrt{\sigma Z + \hat{p}^2 \gamma} = \hat{p}\sqrt{\gamma - \gamma_0}$, where $\gamma_0 := -\sigma Z/\hat{p}^2$ shows the random time the $k$-core first appears.

It is possible that our methods can be developed further to give more precise results for the random graph process. In particular, it would be interesting to find the magnitude of the variations in the size of the first nonempty $k$-core. Are they of the order $n^{1/2}$?

It would also be interesting to understand the rapid growth of the $k$-core in the beginning. Presumably, this happens because there are many rather big subgraphs that have very few vertices with degree less than $k$, and these are quickly joined to and absorbed by the $k$-core. It would be interesting to understand the structure of these subgraphs.

REMARK 1.7. It is well known that $c_2 = 1$ and that when $\lambda > c_2 = 1$, there exists whp a unique giant component in the random graph. The 2-core contains some cycles outside the giant component, but these are too few to influence the asymptotics of its size; thus, from the point of view of this paper, it does not matter whether we study the 2-core of the graph or just the 2-core of the giant component. It was shown by Pittel [26] that, for both $G(n, p)$ and $G(n, m)$, $v(\text{Core}_2)/n \xrightarrow{\text{p}} (1-T)(1-T/\lambda)$, where $T < 1$ satisfies $Te^{-T} = \lambda e^{-\lambda}$. He also conjectured asymptotic normality of $v(\text{Core}_2)$, with an asymptotic variance of the order $n$, which we prove in the present paper.

REMARK 1.8. We can also obtain similar results on the number of vertices of given degree in the core, as has been done by a different method by Cain and Wormald [5]. Indeed, already the simple version of our argument in [17] shows easily that in the case treated in Theorem 1.2 here, if $j \geq k$ and $n_{kj}$ is the number of vertices of degree $j$ in $\mathrm{Core}_k$, then $n_{kj}/n \xrightarrow{\mathrm{P}} \pi_j(\hat{\mu})$, as found by Cain and Wormald [5]. Moreover, it follows from the arguments in this paper that $n^{-1/2}(n_{kj} - \pi_j(\hat{\mu}_n)n)$ has a normal limit, jointly for all $j \geq k$. There is also a similar result in the situation of Theorem 1.3. We omit the details.

REMARK 1.9.   In this paper we treat $G(n,p)$ and $G(n,m)$ together, since all arguments work for both random graphs with no or very minor changes. It is also possible to do the proofs for one version only and then derive the results for the other. To go from results for $G(n,m)$ to $G(n,p)$, we simply condition on the number of edges as in [26]; note that the asymptotic variance will be larger for $G(n,p)$ than for $G(n,m)$ since the mean will shift with the number of edges. It is also possible to go in the opposite direction since the $k$-core grows monotonously with the number of edges; see [15].

**2. Preliminaries.**   It will be convenient to work with *multigraphs*, that is, to allow multiple edges and loops. In particular, we shall use the following type of random multigraph.

Let $n \in \mathbb{N}$ and let $(d_i)_1^n$ be a sequence of nonnegative integers such that $\sum_{i=1}^n d_i$ is even. We define a *random multigraph with given degree sequence* $(d_i)_1^n$, denoted by $G^*(n,(d_i)_1^n)$, by the configuration model (see, e.g., [4]): take a set of $d_i$ *half-edges* for each vertex $i$ (with these sets disjoint), and join the half-edges into edges by a partition of the set of all half-edges into pairs; such a partition of the half-edges is known as a *configuration*, and we choose the configuration at random, uniformly over all possible configurations. Note that conditioned on the multigraph being a simple graph, we obtain a uniformly distributed random graph with the given degree sequence, which we denote by $G(n,(d_i)_1^n)$.

Let $u_r$ denote the number of vertices of degree $r$. We further write $2m := \sum_{i=1}^n d_i = \sum_r r u_r$, so that $m$ is the number of edges in the multigraph $G^*(n,(d_i)_1^n)$.

We will let $n \to \infty$, and assume that we for each $n$ are given $(d_i)_1^n = (d_i^{(n)})_1^n$ satisfying the following regularity conditions; cf. Molloy and Reed [24]. [Our condition (ii) can probably be relaxed, but it will be convenient to work with. Note that it implies $\max_i d_i = o(\log n)$.]

CONDITION 2.1.   *For each $n$, $(d_i)_1^n = (d_i^{(n)})_1^n$ is a sequence of nonnegative integers such that $\sum_{i=1}^n d_i$ is even and, for some probability distribution $(p_r)_{r=0}^\infty$ independent of $n$, with $u_r = u_r(n)$ as defined above:*

(i) $u_r/n := \#\{i : d_i = r\}/n \to p_r$ *for every $r \geq 0$ as $n \to \infty$;*
(ii) *for every $A > 1$, we have $\sum_r u_r A^r = \sum_{i=1}^n A^{d_i} = O(n)$.*

*We further assume $p_0 < 1$.*

Let $D_n$ be a random variable with the distribution $\mathbb{P}(D_n = r) = u_r/n = u_r(n)/n$; this is the distribution of the vertex degree of a random vertex in $G^*(n, (d_i)_1^n)$. Further, let $D$ be a random variable with the distribution $\mathbb{P}(D = r) = p_r$, and let $\lambda := \mathbb{E} D = \sum_r r p_r$. Then Condition 2.1 can be rewritten as

$$D_n \xrightarrow{\mathrm{d}} D \qquad \text{as } n \to \infty, \tag{2.1}$$

$$\mathbb{E} A^{D_n} = O(1) \qquad \text{for each } A > 1, \tag{2.2}$$

and $\lambda > 0$. Note that (2.2) implies uniform integrability of $D_n$, so

$$\frac{2m}{n} = \frac{\sum_{i=1}^n d_i}{n} = \mathbb{E} D_n \to \mathbb{E} D = \lambda = \sum_r r p_r; \tag{2.3}$$

in particular, $\lambda < \infty$. Similarly, all higher moments converge. (It follows that our assumption Condition 2.1 is stronger than the corresponding assumptions in Molloy and Reed [24, 25] and Janson and Luczak [17].)

REMARK 2.2. The excluded case $p_0 = 1$ (where thus $\lambda = 0$) is rather degenerate; some of the results below hold trivially, but others may fail. In this case, most vertices are isolated. By removing all such vertices, we obtain a smaller random graph of the same type, and our results can be applied under suitable conditions.

We use the notation $\mathbf{1}[\mathcal{E}]$ for the indicator of an event $\mathcal{E}$; this is thus 1 if $\mathcal{E}$ holds and 0 otherwise. We also write $\mathbb{E}(X; \mathcal{E})$ for $\mathbb{E}(X \cdot \mathbf{1}[\mathcal{E}])$. We use $(X \mid \mathcal{E})$ to denote a random variable $X$ conditioned on an event $\mathcal{E}$; thus, $\mathbb{E}(X \mid \mathcal{E}) = \mathbb{E}(X; \mathcal{E})/\mathbb{P}(\mathcal{E})$.

We shall consider random thinnings of the vertex degrees. In general, if $X$ is a nonnegative integer valued random variable and $0 \le p \le 1$, we let $X_p$ be the thinning of $X$ obtained by taking $X$ points and then randomly and independently keeping each of them with probability $p$. For integers $l \ge 0$ and $0 \le r \le l$, let $\beta_{lr}$ denote the binomial probabilities

$$\beta_{lr}(p) := \mathbb{P}(\mathrm{Bi}(l,p) = r) = \binom{l}{r} p^r (1-p)^{l-r}.$$

Then we have

$$\mathbb{P}(X_p = r) = \sum_{l=r}^{\infty} \mathbb{P}(X = l) \beta_{lr}(p).$$

We further define the functions, for $0 \leq p \leq 1$, $j = 0, 1, \ldots$, and a given $k \geq 2$,

$$(2.4) \qquad q_{jX}(p) := \mathbb{P}(X_p \geq j) = \sum_{r=j}^{\infty} \sum_{l=r}^{\infty} \mathbb{P}(X = l) \beta_{lr}(p),$$

$$(2.5) \qquad b_X(p) := q_{kX}(p) = \sum_{r=k}^{\infty} \sum_{l=r}^{\infty} \mathbb{P}(X = l) \beta_{lr}(p),$$

$$(2.6) \qquad h_X(p) := \mathbb{E}(X_p; X_p \geq k) = \sum_{r=k}^{\infty} \sum_{l=r}^{\infty} r \mathbb{P}(X = l) \beta_{lr}(p),$$

and, provided $\mathbb{E} X < \infty$,

$$(2.7) \qquad l_X(p) := (\mathbb{E} X) p^2 - h_X(p).$$

These functions for $X = D_n$ and $X = D$ will play an important role in the sequel, and we use the abbreviated notation $q_{jn} = q_{jD_n}$, $q_j = q_{jD}$, $b_n = b_{D_n}$, $b = b_{jD}$, etc. They are thus given by (2.4)–(2.7) with $\mathbb{P}(X = l)$ replaced by $u_l/n$ and $p_l$, and $\mathbb{E} X$ replaced by $2m/n$ and $\lambda$, respectively.

LEMMA 2.3. *As $n \to \infty$, $b_n(p) \to b(p)$, $h_n(p) \to h(p)$ and $l_n(p) \to l(p)$, together with all derivatives, uniformly on $[0, 1]$.*

PROOF. Consider, for example,

$$h_n(p) = \sum_{l=k}^{\infty} \sum_{r=k}^{l} r \beta_{lr}(p) \frac{u_l}{n} = \mathbb{E} f(D_n; p),$$

where $f(l; p) := \sum_{r=k}^{l} r \beta_{lr}(p)$. Then $f(D_n; p) \xrightarrow{d} f(D; p)$ by (2.1). Moreover, $0 \leq f(l; p) \leq l$, so $f(D_n; p) \leq D_n$ is uniformly integrable by (2.2). Hence,

$$h_n(p) = \mathbb{E} f(D_n; p) \to \mathbb{E} f(D; p) = h(p)$$

for each $p$.

Next, an elementary calculation yields

$$(2.8) \qquad \frac{d}{dp} \beta_{lr}(p) = l \beta_{l-1, r-1}(p) - l \beta_{l-1, r}(p).$$

Thus, for every $j \geq 0$, using a simple induction and the fact that $0 \leq \beta_{lr}(p) \leq 1$, we obtain $|\frac{d^j}{dp^j} \beta_{lr}(p)| \leq (2l)^j$. It follows that $|\frac{d^j}{dp^j} f(l, p)| \leq (2l)^j \sum_{r=k}^{l} r \leq (2l)^{j+2}$. Hence,

$$\left| \frac{d^j}{dp^j} f(D_n, p) \right| \leq (2 D_n)^{j+2},$$

which by (2.2) is uniformly integrable, and thus,

$$\frac{d^j}{dp^j}h_n(p) = \mathbb{E}\frac{d^j}{dp^j}f(D_n,p) \to \mathbb{E}\frac{d^j}{dp^j}f(D,p) = \frac{d^j}{dp^j}h(p).$$

Moreover, these derivatives are all uniformly bounded, and uniform convergence for $p \in [0,1]$ follows easily (e.g., using the Arzela–Ascoli theorem). This proves the result for $h$. The result for $l$ follows, and the result for $b$ is proved the same way. □

We note also that by (2.8) (or a combinatorial argument),

$$(2.9) \qquad q'_j(p) = \sum_{l=j}^{\infty} p_l l \beta_{l-1,j-1}(p).$$

It follows that $q'_j(p) > 0$ for all $p \in (0,1)$ unless $q_j$ vanishes identically; hence, $b'(p) > 0$ and $h'(p) > 0$ for $0 < p < 1$ unless $h(p) = 0$ for all $p$.

Finally, let us consider the Poisson case, which is central for us, and establish connections with the functions defined in (1.2).

LEMMA 2.4. *If $X \sim \text{Po}(\lambda)$, then*

$$(2.10) \qquad q_{jX}(p) = \mathbb{P}(X_p \geq j) = \psi_j(\lambda p),$$

$$(2.11) \qquad b_X(p) = q_{kX}(p) = \psi_k(\lambda p),$$

$$(2.12) \qquad h_X(p) = \sum_{j=k}^{\infty} j\, \mathbb{P}(X_p = j) = \sum_{j=k}^{\infty} \frac{(\lambda p)^j}{(j-1)!}e^{-\lambda p} = \lambda p \psi_{k-1}(\lambda p)$$

*and*

$$(2.13) \qquad l_X(p) = \lambda p^2 - h_X(p) = \lambda p(p - \psi_{k-1}(\lambda p)).$$

PROOF. The thinning $X_p \sim \text{Po}(\lambda p)$; thus, (2.4)–(2.7) yield the result. □

**3. Finding the core.** The $k$-core of an arbitrary finite graph or multigraph can be found by removing edges where one endpoint has degree $< k$, until no such edges remain, and finally removing all isolated vertices. The order of removal does not matter, and we choose the edges to be deleted at random as follows.

Regard each edge as consisting of two *half-edges*, each half-edge having one endpoint. Say that a vertex is *light* if its degree is $< k$, and *heavy* if its degree is $\geq k$. Similarly, say that a half-edge is light or heavy when its endpoint is. As long as there are some light half-edges, choose one such half-edge uniformly at random and remove the edge it belongs to. (Note that this may change the other endpoint from heavy to light, and thus create new light half-edges.) When there are no light half-edges left, we stop. Then all light vertices are isolated; the heavy vertices and the remaining edges form the $k$-core of the original graph.

We apply this algorithm to a random multigraph $G^*(n,(d_i)_1^n)$ with given degree sequence $(d_i)_1^n$. We use the configuration model as described in Section 2, and combine the core-finding algorithm with the generation of the random configuration by revealing the pairs in the configuration, and thus the edges in $G^*(n,(d_i)_1^n)$, only when the edges are removed by the algorithm. (The remaining pairs are revealed when the algorithm stops.) It is easily seen (see [17] for details) that if we observe only the vertex degrees in the resulting multigraph process, they can be described by the following process of colored (white or red) balls (representing half-edges) in $n$ bins (representing vertices), if we consider the white balls only. (The careful reader will notice that the version here differs slightly from [17], but it is obviously equivalent.)

Begin with $d_i$ balls in bin $i$, $i=1,\ldots,n$; there are thus $2m$ balls in total. Initially, all balls are white. Say that a bin is *light* (at a given time) if it contains $< k$ white balls, and *heavy* if it contains $\geq k$. Similarly, say that a ball is light or heavy when the bin it belongs to is.

We start by coloring a random light ball red. The process then runs in continuous time such that each ball has an exponentially distributed random lifetime with mean 1, independent of the other balls. This means that balls die and are removed at rate 1, independently of one another. Further, when a white ball dies, a randomly chosen light white ball is colored red, provided that there is some such ball; we stop when a white ball dies and there is no light white ball left. (The interpretation in terms of the graph is that a white ball that is colored red and the next white ball that dies represent two half-edges that are joined to form an edge, and this edge is deleted by the core-finding process.)

We let $L(t)$, $H(t)$ and $B(t)$ denote the numbers of light white balls, heavy balls (which always are white) and heavy bins at time $t$, respectively, and let $\tau$ be the time the process stops. ($B$ is denoted $H_1$ in [17].) There are no white light balls left at $\tau$, which would give $L(\tau)=0$, but the last deletion and recoloring step is incomplete, so we rather define $L(\tau):=-1$, pretending that we did recolor a (nonexisting) ball at $\tau$ too. Moreover, the heavy balls left at $\tau$ are exactly the half-edges in the $k$-core. Hence, the number of edges in the $k$-core is $\frac{1}{2}H(\tau)$, while the number of vertices is $B(\tau)$.

We now consider a sequence $G^*(n,(d_i)_1^n)$, $n\geq 1$, with $(d_i)_1^n=(d_i^{(n)})_1^n$ satisfying Condition 2.1; we use the notation $L_n(t)$, $H_n(t)$, $B_n(t)$ and $\tau_n$. We will change variables and define, for $t\geq 0$,

(3.1) $$\check{b}_n(t):=b_n(e^{-t}),$$

(3.2) $$\check{h}_n(t):=h_n(e^{-t}),$$

(3.3) $$\check{l}_n(t) := l_n(e^{-t}) = \frac{2m}{n}e^{-2t} - \check{h}_n(t)$$

and similarly for $\check{b}(t)$, $\check{h}(t)$, $\check{l}(t)$; note that $\check{b}_n \to \check{b}$, $\check{h}_n \to \check{h}$, and $\check{l}_n \to \check{l}$ on $[0, \infty)$, again uniformly and with all derivatives.

We will show that the processes $B_n(t)$, $H_n(t)$ and $L_n(t)$, suitably normalized, converge jointly to certain Gaussian processes. We have defined the processes on $[0, \tau_n]$ only. It is possible, and sometimes useful, to extend the processes to all $t \geq 0$ by changing the rules after $\tau_n$, but we find it simpler to state the result for the stopped processes, even though this makes the next theorem a bit technical. (See also Remark 5.1.) Recall that, if $\tau$ is a stopping time and $X(t)$ is a stochastic process, then the process stopped at $\tau$ is the process

(3.4) $$X^\tau(t) := X(t \wedge \tau), \qquad 0 \leq t < \infty.$$

[We use the notation $x \wedge y := \min(x, y)$.] Note that $X$ only has to be defined on $[0, \tau]$.

For reasons that will become clear in Section 7, we consider also a modification of our processes where we condition the configuration on containing a given (small) set of pairs of half-edges, which means that the multigraph $G^*(n, (d_i)_1^n)$ contains a given set of edges; we say that these half-edges and edges are *golden*. We run the core-finding process as above, with the modification that we never remove a golden edge. (The final result thus contains all golden edges, so it may be larger than the $k$-core.) Translated into our balls-and-bins process, this means that we now also have some golden balls from the beginning; golden balls are always heavy, do not die and are never recolored. A white ball in the same bin as a golden one is, as before, light if the total number of balls in the urn is $< k$, and heavy otherwise, but we say that any bin containing a golden ball is heavy.

We let $D[0, \infty)$ be the standard space of right-continuous functions with left limits on $[0, \infty)$ equipped with the Skorohod topology; see, for example, [11] or [19], Appendix A.2. In particular, if $f$ is continuous, then $f_n \to f$ in $D[0, \infty)$ if and only if $f_n \to f$ uniformly on every compact subinterval.

THEOREM 3.1. *Assume Condition* 2.1 *and use the notation above. Let* $\tau'_n \leq \tau_n$ *be a stopping time such that* $\tau'_n \xrightarrow{\mathrm{P}} t_0$ *for some* $t_0 \geq 0$. *Then, jointly in* $D[0, \infty)$,

(3.5) $$n^{-1/2}(B_n(t \wedge \tau'_n) - n\check{b}_n(t \wedge \tau'_n)) \xrightarrow{\mathrm{d}} Z_B(t \wedge t_0),$$

(3.6) $$n^{-1/2}(H_n(t \wedge \tau'_n) - n\check{h}_n(t \wedge \tau'_n)) \xrightarrow{\mathrm{d}} Z_H(t \wedge t_0),$$

(3.7) $$n^{-1/2}(L_n(t \wedge \tau'_n) - n\check{l}_n(t \wedge \tau'_n)) \xrightarrow{\mathrm{d}} Z_L(t \wedge t_0),$$

where $Z_B$, $Z_H$ and $Z_L$ are continuous Gaussian processes on $[0, t_0]$ with mean 0 and covariances that satisfy, for $0 \leq t \leq t_0$ and $\nu, \varkappa \in \{B, H, L\}$,

$$(3.8) \qquad \operatorname{Cov}(Z_\nu(t), Z_\varkappa(t)) = \sigma_{\nu\varkappa}(e^{-t}),$$

where $\sigma_{\nu\varkappa}(x) = \sigma_{\nu\varkappa}(x; (p_r)_0^\infty)$ are given by (5.37)–(5.47).

If further $t_0 > 0$ and $\mathbb{P}(D > k) > 0$, then the covariance matrix of $(Z_B(t), Z_H(t), Z_L(t))$ is nonsingular for any fixed $t$ with $0 < t \leq t_0$.

The same results hold if, for each $n$, we select a set of $O(1)$ pairs of golden half-edges.

The proof of Theorem 3.1 is given in Section 5.

In order to specify the distribution of the Gaussian processes $Z_B$, $Z_H$ and $Z_L$ completely, we also need their covariances for a pair of distinct times. These too can be determined from the proof but, for simplicity (and since we do not need them for our further results), we give explicit formulas only in the single time case.

REMARK 3.2. In many cases, for example, Theorem 3.4 below, $\tau_n \overset{\mathrm{P}}{\longrightarrow} t_0$ for some $t_0$, and we may, and should, take $\tau'_n := \tau_n$. The reason for introducing $\tau'_n$ is that there are some interesting cases, in particular, in Theorem 3.5 below, where $\tau_n$ does not converge in probability to a constant; we then stop at a possibly earlier time $\tau'_n$ in order to obtain simple results (and proofs). Note that it was shown in [17] (without golden balls) that we have

$$(3.9) \qquad \sup_{t \leq \tau_n} |L_n(t)/n - \check{l}(t)| \overset{\mathrm{P}}{\longrightarrow} 0;$$

this is less precise but valid up to $\tau_n$. It is easily seen that (3.9) holds also if we allow $O(1)$ golden balls. [This is true since, for instance, (3.9) follows easily from (5.4), (5.12) and (5.15) in the proof below.]

REMARK 3.3. The extra condition $\mathbb{P}(D > k) > 0$ is equivalent to $p_j > 0$ for some $j \geq k + 1$. If this fails, then we have the following three exceptional cases, which also follow from the proof in Section 5 (we omit the details):

(i) If $\mathbb{P}(D > k) = 0$ but $\mathbb{P}(D = k) > 0$ ($p_j = 0$ for $j \geq k+1$ but $p_k > 0$), then $Z_H = kZ_B$ and $(Z_B(t), Z_L(t))$ has a nonsingular Gaussian distribution for $0 < t \leq t_0$, except when $k = 2$ and also $p_1 = 0$ ($D \in \{0, 2\}$ a.s.).

(ii) If $k = 2$ and $D \in \{0, 2\}$ a.s. ($p_j = 0$ for $j \neq 0, 2$), then $Z_H = 2Z_B$ and $Z_L = 0$; further, $Z_B(t)$ has a nonsingular Gaussian distribution for $0 < t \leq t_0$.

(iii) If $\mathbb{P}(D \geq k) = 0$ ($p_j = 0$ for $j \geq k$), then $Z_B = Z_H = 0$, while $Z_L(t)$ has a nonsingular Gaussian distribution for $0 < t \leq t_0$.

From Theorem 3.1 we deduce the following two results for $G(n,(d_i)_1^n)$, which are analogues of Theorems 1.2 and 1.3 for $G(n,p)$ and $G(n,m)$. Again, these two results show the joint asymptotic normality of the number of vertices and edges in the $k$-core above the threshold, as well as a more complicated limit law at the threshold. The proofs are given in Sections 6 and 7.

THEOREM 3.4. *Let $k \geq 2$. Assume Condition 2.1 and use the notation above. Let $\mathrm{Core}_k$ be the $k$-core of $G^*(n,(d_i)_1^n)$. Let $\hat{p}$ be the largest $p \leq 1$ such that $l(p) = 0$, and suppose that $0 < \hat{p} < 1$. Suppose further that $l'(\hat{p}) > 0$. Let $\alpha := l'(\hat{p}) = 2\lambda\hat{p} - h'(\hat{p}) > 0$, let $\hat{t} = -\ln\hat{p}$ and let $Z_L$, $Z_H$ and $Z_B$ be as in Theorem 3.1. Then*

$$\begin{aligned}
(3.10) \quad & n^{-1/2}(v(\mathrm{Core}_k) - b_n(\hat{p}_n)n, e(\mathrm{Core}_k) - \tfrac{1}{2}h_n(\hat{p}_n)n) \\
& \xrightarrow{\mathrm{d}} (Z_B(\hat{t}) - \alpha^{-1}b'(\hat{p})Z_L(\hat{t}), \tfrac{1}{2}Z_H(\hat{t}) - \tfrac{1}{2}\alpha^{-1}h'(\hat{p})Z_L(\hat{t})),
\end{aligned}$$

*where $\hat{p}_n$ is the largest $p \leq 1$ such that $l_n(p) = 0$; further, $\hat{p}_n \to \hat{p}$. The limit distribution is a Gaussian random vector with a nonsingular covariance matrix.*

*The same result holds for $G(n,(d_i)_1^n)$.*

The covariance matrix is easily calculated explicitly, but the formulas are complicated and we refer the reader to (3.8) and (5.37)–(5.47).

THEOREM 3.5. *Let $k \geq 3$. Assume Condition 2.1 and use the notation above. Let $\mathrm{Core}_k$ be the $k$-core of $G^*(n,(d_i)_1^n)$. Suppose that $\min_{0 \leq p \leq 1} l(p) = 0$. Suppose further that there exists a unique $\hat{p} \in (0,1]$ with $l(\hat{p}) = 0$, that $\hat{p} < 1$, and that $\beta := l''(\hat{p}) > 0$. Then, for every $\delta_0 > 0$ with $\delta_0 < \hat{p}$, at least for $n$ sufficiently large, $l_n$ has a unique minimum point $\bar{p}_n$ in $[\delta_0, 1]$, and $\bar{p}_n \to \hat{p}$ and $l_n(\bar{p}_n) \to 0$ as $n \to \infty$. Let further $\sigma := \sigma_{LL}(\hat{p})^{1/2} > 0$, given by (5.47):*

(i) *If $n^{1/2}l_n(\bar{p}_n) \to +\infty$, then whp $\mathrm{Core}_k$ is empty.*
(ii) *If $n^{1/2}l_n(\bar{p}_n) \to \zeta \in (-\infty, \infty)$, then*

$$(3.11) \qquad \mathbb{P}(\mathrm{Core}_k \neq \varnothing) \to \Phi(-\zeta/\sigma),$$

*where $\Phi$ is the standard normal distribution function, and, with $Z \sim N(0,1)$,*

$$\begin{aligned}
& (n^{-3/4}(v(\mathrm{Core}_k) - b_n(\bar{p}_n)n, e(\mathrm{Core}_k) - \tfrac{1}{2}h_n(\bar{p}_n)n)|\mathrm{Core}_k \neq \varnothing) \\
& \to ((2/\beta)^{1/2}\sqrt{\sigma Z - \zeta}(b'(\hat{p}), \tfrac{1}{2}h'(\hat{p})) \mid Z > \zeta/\sigma).
\end{aligned}$$

(iii) *If $n^{1/2}l_n(\bar{p}_n) \to -\infty$, then whp $\mathrm{Core}_k$ is nonempty. Moreover, there is, at least for large $n$, a unique $\hat{p}_n$ with $l_n(\hat{p}_n) = 0$ and $\bar{p}_n < \hat{p}_n < 1$, and*

$$
\begin{aligned}
(3.12) \quad &|l_n(\bar{p}_n)|^{1/2} n^{-1/2}(v(\mathrm{Core}_k) - b_n(\hat{p}_n)n, e(\mathrm{Core}_k) - \tfrac{1}{2}h_n(\hat{p}_n)n) \\
&\xrightarrow{\mathrm{d}} (b'(\hat{p})Z', \tfrac{1}{2}h'(\hat{p})Z'),
\end{aligned}
$$

*where $Z' \sim N(0, \sigma^2/(2\beta))$.*

*The same results hold for $G(n, (d_i)_1^n)$.*

*If further $l'_n(\hat{p})^2 = o(l_n(\hat{p}))$, then in the above results $l_n(\bar{p}_n)$ may be replaced by $l_n(\hat{p})$. Moreover, if $l'_n(\hat{p}) = o(n^{-1/4})$, then $b_n(\bar{p}_n)$ and $h_n(\bar{p}_n)$ may be replaced by $b_n(\hat{p})$ and $h_n(\hat{p})$ in (ii).*

**4. Some martingale theory.** Our proof is based on martingale theory, in particular, a martingale limit theorem by Jacod and Shiryaev [11]. We will use the *quadratic variation* $[X, X]_t$ of a martingale $X$ defined on $[0, \infty)$, and its bilinear extension $[X, Y]_t$ to two martingales $X$ and $Y$. For a general definition see, for example, [11]; for us it will suffice to know that, if $X$ and $Y$ are martingales of pathwise finite variation, then

$$
(4.1) \qquad [X, Y]_t = \sum_{0 < s \leq t} \Delta X(s) \Delta Y(s),
$$

where $\Delta X(s) := X(s) - X(s-)$ is the jump of $X$ at $s$ and, similarly, $\Delta Y(s) := Y(s) - Y(s-)$. The sum in (4.1) is formally uncountable, but in reality countable since there is only a countable number of jumps; in the applications below, the sum will be finite. Note that $[X - X_0, Y - Y_0] = [X, Y]$. (There is some disagreement in the literature concerning the definition of $[X, Y]$ in the case where $X(0)Y(0) \neq 0$; we have chosen the version in [11], with $[X, Y]_0 = 0$.)

For vector-valued martingales $X = (X_i)_{i=1}^m$ and $Y = (Y_j)_{j=1}^n$, we define the *square bracket* $[X, Y]$ to be the $m \times n$ matrix $([X_i, Y_j])_{i,j}$.

A real-valued martingale $X(s)$ on $[0, t]$ is an $L^2$-martingale if and only if $\mathbb{E}[X, X]_t < \infty$ and $\mathbb{E}|X(0)|^2 < \infty$, and then

$$
(4.2) \qquad \mathbb{E}|X(t)|^2 = \mathbb{E}[X, X]_t + \mathbb{E}|X(0)|^2.
$$

We will use the following general result based on [11]; see [14], Proposition 9.1, for a detailed proof. (See also [12] and [13] for similar versions.)

PROPOSITION 4.1. *Assume that for each $n$, $M_n(t) = (M_{ni}(t))_{i=1}^q$ is a $q$-dimensional martingale on $[0, \infty)$ with $M_n(0) = 0$, and that $\Sigma(t)$, $t \geq 0$, is a (nonrandom) continuous matrix-valued function such that, for every fixed $t \geq 0$,*

$$
(4.3) \qquad [M_n, M_n]_t \xrightarrow{\mathrm{p}} \Sigma(t) \qquad \text{as } n \to \infty,
$$

$$
(4.4) \qquad \sup_n \mathbb{E}[M_{ni}, M_{ni}]_t < \infty, \qquad i = 1, \ldots, q.
$$

Then $M_n \xrightarrow{d} M$ as $n \to \infty$, in $D[0, \infty)$, where $M$ is a continuous $q$-dimensional Gaussian martingale with $\mathbb{E} M(t) = 0$ and covariances

$$\mathbb{E} M(t) M'(u) = \Sigma(t), \qquad 0 \leq t \leq u < \infty.$$

REMARK 4.2. By (4.2), (4.4) is equivalent to $\sup_n \mathbb{E} |M_n(t)|^2 < \infty$, the form used in, for example, [14].

The proposition thus yields joint convergence of the processes $M_{ni}$, $1 \leq i \leq q$. This extends immediately to infinitely many processes (formally the case $q = \infty$) by considering finite subsets, since, by definition, an infinite family of random variables (or processes) converge jointly if every finite subfamily does.

We will apply Proposition 4.1 to stopped processes. We say that a stochastic process $X$ is a martingale on $[0, \tau]$, where $\tau$ is a stopping time, if the stopped process $X^\tau$ defined in (3.4) is a martingale on $[0, \infty)$.

Throughout our proofs we shall use $\mathcal{M}$ or $\mathcal{M}'$ to denote an unspecified martingale.

**5. Proof of Theorem 3.1.** We assume for simplicity that there are no golden balls, and leave to the reader the very minor modifications in the case with $O(1)$ golden balls. (When defining $U_{rn}$ and $V_{rn}$ below, we count only those bins that do not contain any golden balls.)

Consider first $W_n(t) := L_n(t) + H_n(t)$, the number of white balls. By construction, $W_n(0) = 2m - 1$. Moreover, $W_n(t)$ decreases by 2 each time a white ball dies. Since each ball dies with rate 1, it follows that on $[0, \tau_n]$, writing in differential form,

(5.1) $$dW_n(t) = -2W_n(t)\, dt + d\mathcal{M}(t).$$

In other words, $W_n(t) + \int_0^t 2W_n(s)\, ds$ is a martingale on $[0, \tau_n]$. In this section we will for simplicity usually omit the qualification "on $[0, \tau_n]$," but it is tacitly assumed that we consider $t \leq \tau_n$ only, unless stated otherwise.

REMARK 5.1. It is possible to continue the process $W_n$ after $\tau_n$, so that (5.1) still holds, by merging all the bins into one and using the straightforward rule that balls die with rate 1 and that, whenever a ball dies, a second ball is removed. This works until there is only one ball left; we may also change $W_n$ by $\pm 1$ as in [17] and have $W_n(t)$ defined for all $t \geq 0$, still with (5.1). Similarly, we may extend the $U_{rn}$ and all the other processes derived from it below to all $t \geq 0$; the extended processes are defined by ignoring the colors of remaining balls, so that we simply have bins with balls that die independently of one another at rate 1. However, we have been unable to find a good way to extend both $W_n$ and the other processes together. We therefore usually consider the processes up to $\tau_n$ only.

By (5.1) and Itô's lemma,
$$d(e^{2t}W_n(t)) = e^{2t}\,dW_n(t) + 2e^{2t}W_n(t)\,dt = e^{2t}\,d\mathcal{M}(t),$$
another martingale differential, that is, $\widehat{W}_n(t) := e^{2t}W_n(t)$ is a martingale. (Note that, for each $n$, $W_n(t)$ is bounded by $2m(n)$. Thus, $\widehat{W}_n(t)$ is bounded on every finite interval. Hence, $\widehat{W}_n(t \wedge \tau_n)$ is a square integrable martingale on every interval $[0,T]$, and not just a local martingale. The same applies to the other martingales in this section.)

Since distinct balls a.s. die at distinct times, all jumps in $W_n(t)$ equal $-2$, and the quadratic variation of $\widehat{W}_n$ is given by, and, using integration by parts, the quadratic variation

$$[\widehat{W}_n, \widehat{W}_n]_t = \sum_{0<s\leq t} |\Delta \widehat{W}_n(s)|^2 = \sum_{0<s\leq t} |e^{2s}\Delta W_n(s)|^2 = \int_0^t 2e^{4s}\,d(-W_n(s))$$
(5.2)
$$= -2e^{4t}W_n(t) + 2W_n(0) + \int_0^t 8e^{4s}W_n(s)\,ds.$$

In particular,

$$(5.3) \qquad [\widehat{W}_n, \widehat{W}_n]_t \leq 2e^{4t} \int_0^t d(-W_n(s)) \leq 2e^{4t}W_n(0) \leq 4me^{4t}.$$

Let $W_n^*(t) := n^{-1}\widehat{W}_n(t \wedge \tau_n)$. Then $W_n^*$ is a martingale on $[0,\infty)$ and, for every $T < \infty$,
$$[W_n^*, W_n^*]_T = n^{-2}[\widehat{W}_n, \widehat{W}_n]_{T \wedge \tau_n} \leq 4e^{4T}m/n^2 \to 0$$
as $n \to \infty$, by (2.3). Moreover, $0 \leq W_n^*(t) \leq 2m/n = O(1)$. We can thus apply Proposition 4.1 to $W_n^*(t) - W_n^*(0)$, with $q=1$ and $\Sigma(t) = 0$. In this case, the limit $M$ satisfies $\mathrm{Var}(M(t)) = \Sigma(t) = 0$, and thus, $M(t) = 0$ a.s., for every $t$. Consequently, we have shown $W_n^*(t) - W_n^*(0) \xrightarrow{\mathrm{P}} 0$ in $D[0,\infty)$, and thus, $W_n^*(t) - W_n^*(0) \xrightarrow{\mathrm{P}} 0$ uniformly on $[0,T]$ for every $T < \infty$. Equivalently,

$$(5.4) \qquad n^{-1} \sup_{t \leq T} |W_n(t) - W_n(0)e^{-2t}| \xrightarrow{\mathrm{P}} 0.$$

(This was shown in [17], (5.1) by other methods. Actually, [17], (5.1) is equivalent to (5.4) above with $T = \infty$, which is an easy consequence of (5.4) for finite $T$ because $W_n(t)$ is decreasing and $e^{-2t} \to 0$ as $t \to \infty$.)

Returning to (5.2), we obtain by (5.4), for every $t \in [0, T \wedge \tau_n]$ with $T$ fixed,

$$[\widehat{W}_n, \widehat{W}_n]_t = -2e^{4t}W_n(0)e^{-2t} + 2W_n(0) + \int_0^t 8e^{4s}W_n(0)e^{-2s}\,ds + o_\mathrm{p}(n)$$
$$= W_n(0)(-2e^{2t} + 2 + 4e^{2t} - 4) + o_\mathrm{p}(n)$$
$$= 2W_n(0)(e^{2t} - 1) + o_\mathrm{p}(n) = 2(2m-1)(e^{2t} - 1) + o_\mathrm{p}(n)$$
$$= 2\lambda n(e^{2t} - 1) + o_\mathrm{p}(n).$$

Defining
$$\widetilde{\widehat{W}}_n(t) := n^{-1/2}(\widehat{W}_n(t) - \widehat{W}_n(0)) = n^{-1/2}(e^{2t}W_n(t) - W_n(0)),$$

we have

(5.5) $\quad [\widetilde{\widehat{W}}_n, \widetilde{\widehat{W}}_n]_t = n^{-1}[\widehat{W}_n, \widehat{W}_n]_t = 2\lambda(e^{2t} - 1) + o_{\mathrm{p}}(1), \qquad t \leq \tau_n.$

Let us now stop the processes at $\tau'_n \leq \tau_n$. By (5.5) and the assumption $\tau'_n \xrightarrow{\mathrm{p}} t_0$, the quadratic variation of the stopped process converges in probability to $2\lambda(e^{2(t \wedge t_0)} - 1)$, for any fixed $t \in [0, \infty)$. Moreover, by (5.3),
$$[\widetilde{\widehat{W}}_n, \widetilde{\widehat{W}}_n]_{t \wedge \tau'_n} = n^{-1}[\widehat{W}_n, \widehat{W}_n]_{t \wedge \tau'_n} \leq 4me^{4t}/n,$$

which is bounded for every fixed $t$. Consequently, Proposition 4.1 applies to the stopped process and shows that

(5.6) $\qquad \widetilde{\widehat{W}}_n(t \wedge \tau'_n) \xrightarrow{\mathrm{d}} \widehat{Z}(t \wedge t_0) \qquad \text{in } D[0, \infty),$

where $\widehat{Z}$ is a continuous Gaussian martingale with $\mathbb{E}\,\widehat{Z}(t) = 0$ and covariances

(5.7) $\qquad \mathbb{E}(\widehat{Z}(t)\widehat{Z}(u)) = 2\lambda(e^{2t} - 1), \qquad 0 \leq t \leq u < \infty.$

[We can, if we want to, assume that $\widehat{Z}$ is defined, with this covariance function, for all $t \geq 0$; $\widehat{Z}$ is just a time-change of a standard Brownian motion B: $\widehat{Z}(t) = B(2\lambda(e^{2t} - 1))$.] We define further

(5.8) $\qquad \widetilde{W}_n(t) := n^{-1/2}(W_n(t) - W_n(0)e^{-2t}) = e^{-2t}\widetilde{\widehat{W}}_n(t),$

(5.9) $\qquad Z_W(t) := e^{-2t}\widehat{Z}(t),$

and note that $Z_W$ is a continuous Gaussian process with $Z_W(t) \sim N(0, 2\lambda(e^{-2t} - e^{-4t}))$. Then (5.6) implies

(5.10) $\qquad \widetilde{W}_n(t \wedge \tau'_n) \xrightarrow{\mathrm{d}} Z_W(t \wedge t_0) \qquad \text{in } D[0, \infty).$

Next let us ignore the colors. For $r = 0, 1, \ldots,$ let $U_{rn}(t)$ be the number of bins with exactly $r$ balls at time $t$, and let $V_{rn}(t) := \sum_{s \geq r} U_{sn}(t)$ be the number of bins with at least $r$ balls. Thus, $V_{0n}(t) = \sum_r U_{rn}(t) = n$ and $U_{rn}(0) = u_r(n)$. Note that

(5.11) $\qquad B_n(t) = V_{kn}(t),$

(5.12) $\qquad H_n(t) = \sum_{r=k}^{\infty} rU_{rn}(t) = kV_{kn}(t) + \sum_{k+1}^{\infty} V_{rn}(t).$

Since $V_{rn}$ changes (by $-1$) precisely when a ball dies in a bin with $r$ balls, and there are $rU_{rn}$ such balls,
$$dV_{rn}(t) = -rU_{rn}(t)\,dt + d\mathcal{M}'.$$

Define further $\widehat{V}_{rn}(t) := e^{rt}V_{rn}(t)$. Then

$$\begin{aligned}d\widehat{V}_{rn}(t) &= re^{rt}V_{rn}(t)\,dt + e^{rt}\,dV_{rn}(t) \\ &= re^{rt}V_{rn}(t)\,dt - re^{rt}U_{rn}(t)\,dt + e^{rt}\,d\mathcal{M}' \\ &= re^{rt}V_{r+1,n}(t)\,dt + d\mathcal{M} = re^{-t}\widehat{V}_{r+1,n}(t)\,dt + d\mathcal{M},\end{aligned}$$

where $d\mathcal{M} = e^{rt}\,d\mathcal{M}'$ is another martingale differential. It follows that

$$(5.13) \qquad M_{rn}(t) := \widehat{V}_{rn}(t) - r\int_0^t e^{-s}\widehat{V}_{r+1,n}(s)\,ds$$

is a martingale for every $r \geq 0$. The quadratic variation is

$$(5.14)\begin{aligned}[M_{rn},M_{rn}]_t &= \sum_{0<s\leq t}|\Delta M_{rn}(s)|^2 = \sum_{0<s\leq t}|\Delta \widehat{V}_{rn}(s)|^2 \\ &= \sum_{0<s\leq t}e^{2rs}|\Delta V_{rn}(s)|^2 = \int_0^t e^{2rs}\,d(-V_{rn}(s)).\end{aligned}$$

Define

$$\widetilde{M}_{rn}(t) := n^{-1/2}(M_{rn}(t) - M_{rn}(0));$$

this is a martingale with $\widetilde{M}_{rn}(0) = 0$. Note that

$$M_{rn}(0) = \widehat{V}_{rn}(0) = V_{rn}(0) = \sum_{l \geq r} u_l(n).$$

It follows from [17], Lemma 4.4, that, for every $r \geq 0$, as $n \to \infty$,

$$(5.15) \qquad \sup_{t \geq 0}|V_{rn}(t)/n - q_r(e^{-t})| \xrightarrow{\text{p}} 0,$$

where $q_r(p) := \mathbb{P}(D_p \geq r)$; see Section 2. [This could also be proved similarly to (5.4) above, using (5.23) below.] Using integration by parts twice, we obtain from (5.14) that, if $0 < T < \infty$ is fixed, then, for $t \leq T \wedge \tau_n$,

$$(5.16)\begin{aligned}[\widetilde{M}_{rn},\widetilde{M}_{rn}]_t &= n^{-1}[M_{rn},M_{rn}]_t \\ &= n^{-1}\left(V_{rn}(0) - e^{2rt}V_{rn}(t) + \int_0^t V_{rn}(s)2re^{2rs}\,ds\right) \\ &= q_r(1) - e^{2rt}q_r(e^{-t}) + \int_0^t q_r(e^{-s})2re^{2rs}\,ds + o_{\text{p}}(1) \\ &= \int_0^t e^{2rs}\,d(-q_r(e^{-s})) + o_{\text{p}}(1) = \int_{e^{-t}}^1 p^{-2r}\,dq_r(p) + o_{\text{p}}(1).\end{aligned}$$

Moreover, (5.14) further implies

$$(5.17) \qquad [\widetilde{M}_{rn},\widetilde{M}_{rn}]_t = n^{-1}[M_{rn},M_{rn}]_t \leq n^{-1}e^{2rt}V_{rn}(0) \leq e^{2rt}.$$

Consequently, again by Proposition 4.1 applied to the processes stopped at $\tau'_n$,

$$\widetilde{M}_{rn}(t \wedge \tau'_n) \xrightarrow{\mathrm{d}} Y_r(t \wedge t_0) \qquad \text{in } D[0, \infty), \tag{5.18}$$

where $Y_r$ is a continuous Gaussian martingale with $\mathbb{E} Y_r(t) = 0$ and covariances

$$\mathbb{E}(Y_r(t) Y_r(u)) := \int_{e^{-t}}^{1} p^{-2r} \, dq_r(p), \qquad 0 \le t \le u. \tag{5.19}$$

Furthermore, since $V_{rn}$ and $V_{sn}$, with $s \ne r$, a.s. never change simultaneously, $[\widetilde{M}_{rn}, \widetilde{M}_{sn}] = 0$. Hence, by Proposition 4.1 applied to the vector-valued martingale $(\widetilde{M}_{rn})_{r=0}^R$, for fixed but arbitrary $R \ge 1$, (5.18) holds jointly for all $r \ge 0$, with a diagonal covariance matrix for $(Y_r)_1^\infty$, which implies that the processes $Y_r$ ($r = 0, 1, \ldots$) are all independent.

To deduce results for $V_{rn}$, we have to invert (5.13). We have, using (5.13) repeatedly,

$$\widehat{V}_{rn}(t) = M_{rn}(t) + r \int_0^t e^{-s} \widehat{V}_{r+1,n}(s) \, ds$$

$$= M_{rn}(t) + \int_{s<t} r e^{-s} M_{r+1,n}(s) \, ds$$

$$+ \int_{s_2 < s_1 < t} r e^{-s_1} (r+1) e^{-s_2} \widehat{V}_{r+2,n}(s_2) \, ds_2 \, ds_1,$$

and so on, and it is easily verified by backward induction that

$$\widehat{V}_{rn}(t) = M_{rn}(t)$$
$$+ \sum_{j=r+1}^{\infty} r \binom{j-1}{r} \int_0^t (e^{-s} - e^{-t})^{j-r-1} e^{-s} M_{j,n}(s) \, ds. \tag{5.20}$$

Note that the sum actually is finite for every $n$, since $V_{rn}(t) = 0$ and $M_{rn}(t) = 0$ when $r > n$.

Define, for $t \ge 0$,

$$\hat{v}_{rn}(t) := M_{rn}(0)$$
$$+ \sum_{j=r+1}^{\infty} r \binom{j-1}{r} \int_0^t (e^{-s} - e^{-t})^{j-r-1} e^{-s} M_{j,n}(0) \, ds. \tag{5.21}$$

We claim that $\hat{v}_{rn}(t) = n e^{rt} q_{rn}(e^{-t})$. To see this, it is convenient to temporarily assume that $\widehat{V}_{rn}(t)$ and $M_{rn}(t)$ are defined for all $t \ge 0$; see Remark 5.1. Since then $M_{rn}$ is a martingale on $[0, \infty)$, we have $\mathbb{E} M_{rn}(t) = M_{rn}(0)$, and thus, $\mathbb{E} \widehat{V}_{rn}(t) = \hat{v}_{rn}(t)$. On the other hand, $\mathbb{E} \widehat{V}_{rn}(t) = e^{rt} \mathbb{E} V_{rn}(t)$, and $V_{rn}(t)$ is the number of bins with at least $r$ balls at time $t$ in the process where balls die independently with rate 1 (without stopping). At time 0, the number of balls in a random bin has the distribution of $D_n$, and at time $t$ it thus has the same distribution as the thinned variable $D_{n,e^{-t}}$. Consequently,

$$(5.22) \qquad \hat{v}_{rn}(t) = \mathbb{E}\,\widehat{V}_{rn}(t) = e^{rt}n\,\mathbb{P}(D_{n,e^{-t}} \geq r) = ne^{rt}q_{rn}(e^{-t}),$$

verifying our claim.

Accordingly, for $t \leq \tau_n$, let

$$\widetilde{\widehat{V}}_{rn}(t) := n^{-1/2}(\widehat{V}_{rn}(t) - ne^{rt}q_{rn}(e^{-t}));$$

we find from (5.20)–(5.22) that

$$(5.23) \quad \begin{aligned} \widetilde{\widehat{V}}_{rn}(t) &= \widetilde{M}_{rn}(t) \\ &\quad + \sum_{j=r+1}^{\infty} r\binom{j-1}{r}\int_0^t (e^{-s}-e^{-t})^{j-r-1}e^{-s}\widetilde{M}_{j,n}(s)\,ds. \end{aligned}$$

We can now apply (5.18) (with joint convergence) to any partial sum. Further, from (5.17) we have $[\widetilde{M}_{rn}, \widetilde{M}_{rn}]_t \leq e^{2rt}V_{rn}(0)/n$. By Condition 2.1(ii), for any $A > 1$, there exists $C_A$ such that $V_{rn}(0) \leq A^{-r}\sum_{j=0}^{\infty} u_j A^j \leq C_A A^{-r}n$. Thus, for any $T > 0$, choosing $A = e^{2T+2}$,

$$[\widetilde{M}_{rn}, \widetilde{M}_{rn}]_T \leq e^{2rT}V_{rn}(0)/n \leq C_A e^{-2r}.$$

Hence, using (4.2) and Doob's $L^2$-inequality, we obtain

$$\mathbb{E}\sup_{t \leq T}\widetilde{M}_{rn}^2(t) \leq 4\,\mathbb{E}[\widetilde{M}_{rn}, \widetilde{M}_{rn}]_T \leq C'_T e^{-2r}.$$

Then by the Cauchy–Schwarz inequality,

$$(5.24) \qquad \mathbb{E}\sup_{t \leq T}|\widetilde{M}_{rn}(t)| \leq C''_T e^{-r}.$$

Hence, by (5.18) and Fatou's lemma,

$$(5.25) \qquad \mathbb{E}\sup_{t \leq t_0}|Y_r(t)| \leq Ce^{-r}.$$

Let $R_{Nn}(t)$ be the tail of the sum in (5.23), summing over $j > N$ only. Using (5.24), it is easily seen that, for any fixed $r$ and $T$, $\mathbb{E}\sup_{t \leq T}|R_{Nn}(t)| \to 0$ as $N \to \infty$, uniformly in $n$. Then using the convergence of the partial sums, we may by [2], Theorem 4.2, take the limit (5.18) (in distribution) under the summation sign in (5.23). It follows that

$$(5.26) \qquad \widetilde{\widehat{V}}_{rn}(t \wedge \tau'_n) \xrightarrow{d} \widehat{X}_r(t \wedge t_0),$$

in $D[0,\infty)$ for each $r$, where

$$(5.27) \quad \widehat{X}_r(t) := Y_r(t) + \sum_{j=r+1}^{\infty} r\binom{j-1}{r} \int_0^t (e^{-s} - e^{-t})^{j-r-1} e^{-s} Y_j(s)\, ds.$$

It is an easy consequence of (5.25) that the sum in (5.27) a.s. converges uniformly in $t \leq t_0$, which implies that each $\widehat{X}_r$ is continuous.

Next, define for each $r$

$$(5.28) \qquad \widetilde{V}_{rn}(t) := n^{-1/2}(V_{rn}(t) - nq_{rn}(e^{-t})) = e^{-rt}\widehat{\widetilde{V}}_{rn}(t),$$

$$(5.29) \qquad X_r(t) := e^{-rt}\widehat{X}_r(t).$$

Then, by (5.26),

$$(5.30) \qquad \widetilde{V}_{rn}(t \wedge \tau_n') \xrightarrow{d} X_r(t \wedge t_0).$$

Furthermore, our proof shows also that there is joint convergence for different $r$ in (5.30). Moreover, the same argument as above (with truncation to finite sums, using (5.24) and [2], Theorem 4.2) shows that we can sum over $r \geq k$ to yield, jointly with (5.30) and with a continuous limit on the right-hand side,

$$(5.31) \qquad \sum_{r=k}^{\infty} \widetilde{V}_{rn}(t \wedge \tau_n') \xrightarrow{d} \sum_{r=k}^{\infty} X_r(t \wedge t_0).$$

It now follows from (5.11) and (5.12) that (3.5) and (3.6) hold with $Z_B = X_k$ and $Z_H = kX_k + \sum_{r=k+1}^{\infty} X_r$; see (3.1), (3.2), (2.5) and (2.6).

Finally, let us apply Proposition 4.1 to $\widetilde{M}_{rn}$, for all $r \geq k$, and $\widetilde{\widehat{W}}_n$ together. Each time $V_{rn}$ jumps (necessarily by $-1$) for some $r \geq k$, a white ball dies and thus $W_n$ jumps too (by $-2$). (We restrict ourselves to $r \geq k$ here, and recall that every heavy ball is white.) Thus, the quadratic covariation is

$$[\widetilde{M}_{rn}, \widetilde{\widehat{W}}_n]_t = n^{-1} \sum_{0<s\leq t} \Delta M_{rn}(s) \Delta \widehat{W}_n(s) = n^{-1} \sum_{0<s\leq t} \Delta \widehat{V}_{rn}(s) \Delta \widehat{W}_n(s)$$

$$= n^{-1} \sum_{0<s\leq t} e^{(r+2)s} \Delta V_{rn}(s) \Delta W_n(s) = n^{-1} \int_0^t 2e^{(r+2)s}\, d(-V_{rn}(s)),$$

and [using integration by parts as in (5.16)] we obtain from (5.15) that

$$[\widetilde{M}_{rn}, \widetilde{\widehat{W}}_n]_t = 2\int_0^t e^{(r+2)s}\, d(-q_r(e^{-s})) + o_{\mathrm{p}}(1)$$
$$(5.32)$$
$$= 2\int_{e^{-t}}^1 p^{-r-2}\, dq_r(p) + o_{\mathrm{p}}(1).$$

Hence, Proposition 4.1 implies joint convergence for $\widetilde{M}_{rn}$ ($r \geq k$) and $\widetilde{W}_n$ stopped at $\tau'_n$, and, by (5.23), this also holds for $\widehat{V}_{rn}$. Hence, we can regard the Gaussian processes $Z_W$, $Y_r$ and $X_r$, $r \geq k$, as jointly defined, and jointly Gaussian, and the limits above holding jointly (at least for $r \geq k$). But $L_n(t) = W_n(t) - H_n(t)$ so, combining together (5.11), (5.12), (5.10), (5.30), (5.31), (3.3) and (5.8) yields (3.5)–(3.7), with

$$(5.33) \qquad Z_B := X_k,$$

$$(5.34) \qquad Z_H := kX_k + \sum_{r=k+1}^{\infty} X_r,$$

$$(5.35) \qquad Z_L := Z_W - Z_H = Z_W - kX_k - \sum_{r=k+1}^{\infty} X_r.$$

Next we compute the covariances of the processes $Z_L$, $Z_H$, $Z_B$. As stated above, for simplicity, we shall only consider the case of a single time $t$; we leave the general case of two different times to the reader. For convenience we make the change of variable $t = -\ln x$; thus, $x = e^{-t}$ decreases from 1 to 0. In what follows, we assume that $0 \leq t \leq t_0$ (so that $e^{-t_0} \leq x \leq 1$), and that $r, s \geq k$. Also, for convenience, in (5.42)–(5.47) our results are stated in terms of quantities $\sigma_{WW}(x)$, $\sigma_{rW}(x)$ and $\sigma_{rs}(x)$ defined below.

First, by (5.7),

$$(5.36) \qquad \operatorname{Var}(\widehat{Z}(-\ln x)) = \mathbb{E}(\widehat{Z}(-\ln x)^2) = 2\lambda(x^{-2} - 1)$$

and, thus,

$$(5.37) \quad \sigma_{WW}(x) := \operatorname{Var}(Z_W(-\ln x)) = \operatorname{Var}(x^2 \widehat{Z}(-\ln x)) = 2\lambda(x^2 - x^4).$$

Similarly, by (5.19) and the ensuing comments,

$$\operatorname{Cov}((Y_r(-\ln x), Y_s(-\ln x)) = \delta_{rs} \int_x^1 p^{-2r} \, dq_r(p).$$

Moreover, (5.32) implies, by Proposition 4.1,

$$(5.38) \qquad \operatorname{Cov}(Y_r(-\ln x), \widehat{Z}(-\ln x)) = 2 \int_x^1 p^{-r-2} \, dq_r(p).$$

Let $s = -\ln y$. Recall that $\widehat{Z}$ is a martingale, and so $\operatorname{Cov}(\widehat{Z}(t), Y_j(s)) = \operatorname{Cov}(\widehat{Z}(s), Y_j(s))$. Then from (5.27), for $s \leq t$,

$$\operatorname{Cov}(\widehat{X}_r(-\ln x), \widehat{Z}(-\ln x)) = 2 \int_x^1 p^{-r-2} \, dq_r(p) + \sum_{j=r+1}^{\infty} r \binom{j-1}{r} a_{rj}(x),$$

where we define, with the change of variable $y = e^{-s}$,

$$a_{rj}(x) := \int_{y=x}^{1} (y-x)^{j-r-1} 2 \int_{p=y}^{1} p^{-j-2} \, dq_j(p) \, dy$$

$$= 2 \int_{p=x}^{1} \int_{y=x}^{p} (y-x)^{j-r-1} p^{-j-2} \, dy \, dq_j(p)$$

$$= \frac{2}{j-r} \int_{x}^{1} (p-x)^{j-r} p^{-j-2} \, dq_j(p).$$

Hence,

(5.39)
$$\operatorname{Cov}(\widehat{X}_r(-\ln x), \widehat{Z}(-\ln x))$$
$$= 2 \sum_{j=r}^{\infty} \binom{j-1}{r-1} \int_{x}^{1} (p-x)^{j-r} p^{-j-2} \, dq_j(p)$$

and, recalling (5.9) and (5.29),

(5.40)
$$\sigma_{rW}(x) := \operatorname{Cov}(X_r(-\ln x), Z_W(-\ln x))$$
$$= 2x^{r+2} \sum_{j=r}^{\infty} \binom{j-1}{r-1} \int_{x}^{1} (p-x)^{j-r} p^{-j-2} \, dq_j(p).$$

Similarly, by (5.27) and (5.19),

$$\operatorname{Var}(\widehat{X}_r(-\ln x))$$
$$= \int_{x}^{1} p^{-2r} \, dq_r(p) + \sum_{j=r+1}^{\infty} r^2 \binom{j-1}{r}^2$$
$$\times 2 \iint_{x<y<z<1} (y-x)^{j-r-1} (z-x)^{j-r-1}$$
$$\times \int_{p=z}^{1} p^{-2j} \, dq_j(p) \, dy \, dz$$
$$= \sum_{j=r}^{\infty} \binom{j-1}{r-1}^2 \int_{x}^{1} (p-x)^{2j-2r} p^{-2j} \, dq_j(p)$$

and, omitting the details, for $i \geq 1$,

$$\operatorname{Cov}(\widehat{X}_r(-\ln x), \widehat{X}_{r+i}(-\ln x))$$
$$= \sum_{j=r+i}^{\infty} \binom{j-1}{r-1} \binom{j-1}{r+i-1} \int_{x}^{1} (p-x)^{2j-2r-i} p^{-2j} \, dq_j(p).$$

Hence, by (5.29) again, for $i \geq 0$,

$$\sigma_{r,r+i}(x) := \mathrm{Cov}(X_r(-\ln x), X_{r+i}(-\ln x))$$

(5.41)
$$= x^{2r+i} \sum_{j=r+i}^{\infty} \binom{j-1}{r-1}\binom{j-1}{r+i-1}$$
$$\times \int_x^1 (p-x)^{2j-2r-i} p^{-2j} \, dq_j(p).$$

The covariances of the processes $Z_L$, $Z_H$, $Z_B$ now follow from (5.33)–(5.35), (5.37), (5.40) and (5.41); (3.8) now follows, with

(5.42) $\quad \sigma_{BB}(x) = \sigma_{kk}(x),$

(5.43) $\quad \sigma_{BH}(x) = k\sigma_{kk}(x) + \sum_{r=k+1}^{\infty} \sigma_{kr}(x),$

(5.44) $\quad \sigma_{HH}(x) = k^2 \sigma_{kk}(x) + 2k \sum_{r=k+1}^{\infty} \sigma_{kr}(x) + \sum_{r,s=k+1}^{\infty} \sigma_{rs}(x),$

(5.45) $\quad \sigma_{BL}(x) = \sigma_{kW}(x) - \sigma_{BH}(x),$

(5.46) $\quad \sigma_{HL}(x) = k\sigma_{kW}(x) + \sum_{r=k+1}^{\infty} \sigma_{rW}(x) - \sigma_{HH}(x),$

(5.47) $\quad \sigma_{LL}(x) = \sigma_{WW}(x) - 2k\sigma_{kW}(x) - 2\sum_{r=k+1}^{\infty} \sigma_{rW}(x) + \sigma_{HH}(x).$

Finally, to see that this gives a nonsingular covariance matrix when $\mathbb{P}(D > k) > 0$ and $t = -\ln x \in (0, t_0]$, we argue as follows (without using the explicit formulas above). Assume that the matrix is singular for $t = t_1 > 0$. Then $a_B Z_B(t_1) + a_H Z_H(t_1) + a_L Z_L(t_1) = 0$ a.s. for some constants $a_B, a_H, a_L$, not all zero.

If $a_L = 0$, then (5.33) and (5.34) show that either $X_k = a\sum_{j \geq k+1} X_j$ for some $a$ or $X_{k+1} = -\sum_{j \geq k+2} X_j$ a.s. By (5.29) and (5.27), this means that for either $l = k$ or $l = k+1$, $Y_l$ is a.s. equal to a function of $\{Y_j\}_{j>l}$ that can be written as a sum $\sum_{l+1}^{\infty} \int_0^{t_1} f_j(s) \, dY_j(s)$ of stochastic integrals, for some deterministic functions $f_j(s)$ (depending on $t_1$). However, the processes $Y_j$ are independent, so this is impossible unless $Y_l(t_1) = 0$ a.s., but this contradicts (5.19) and (2.9) since $\mathbb{P}(D \geq k+1) > 0$.

If $a_L \neq 0$, we similarly see that $\widehat{Z}(t_1) = e^{2t_1} Z_W(t_1)$ is a.s. equal to such a sum, now over $j \geq k$, so that for suitable deterministic functions $f_j$,

$$\widehat{Z}(t_1) = \int_0^{t_1} d\widehat{Z}(s) = \sum_{j=k}^{\infty} \int_0^{t_1} f_j(s) \, dY_j(s).$$

(The integrals on the right-hand side are independent, and the sum converges in $L^2$.) In other words, when $t = t_1$,

$$\int_0^t \left( d\widehat{Z}(s) - \sum_{j=k}^{\infty} f_j(s) \, dY_j(s) \right) = 0. \tag{5.48}$$

This stochastic integral is a martingale, and thus, (5.48) holds for all $t \in [0, t_1]$.

Rewrite (5.36)–(5.38) as

$$\operatorname{Var}(\widehat{Z}(t)) = \int_0^t b_{ZZ}(s) \, ds,$$

$$\operatorname{Cov}(Y_j(t), Y_r(t)) = \delta_{jr} \int_0^t b_{jj}(s) \, ds,$$

$$\operatorname{Cov}(\widehat{Z}(t), Y_j(t)) = \int_0^t b_{Zj}(s) \, ds,$$

where

$$b_{ZZ}(t) = 4\lambda e^{2t}, \tag{5.49}$$

$$b_{jj}(t) = e^{(2j-1)t} q'_j(e^{-t}), \tag{5.50}$$

$$b_{Zj}(t) = 2e^{(j+1)t} q'_j(e^{-t}). \tag{5.51}$$

Then, computing the variance of the stochastic integral in (5.48), we find, for $0 \leq t \leq t_1$,

$$\int_0^t \left( b_{ZZ}(s) - 2 \sum_{j=k}^{\infty} f_j(s) b_{Zj}(s) + \sum_{j=k}^{\infty} f_j(s)^2 b_{jj}(s) \right) = 0,$$

and thus, for a.e. $s \in [0, t_1]$,

$$b_{ZZ}(s) = 2 \sum_{j=k}^{\infty} f_j(s) b_{Zj}(s) - \sum_{j=k}^{\infty} f_j(s)^2 b_{jj}(s). \tag{5.52}$$

Now, the arithmetic–geometric inequality shows that

$$2 f_j(s) b_{Zj}(s) \leq f_j(s)^2 b_{jj}(s) + b_{Zj}(s)^2 / b_{jj}(s),$$

so, using formulae (5.49)–(5.51), we deduce from (5.52) that

$$4\lambda e^{2s} = b_{ZZ}(s) \leq \sum_{j=k}^{\infty} b_{Zj}(s)^2 / b_{jj}(s) = 4 \sum_{j=k}^{\infty} e^{3s} q'_j(e^{-s}). \tag{5.53}$$

On the other hand, for any finite sequence $g_k, \ldots, g_N$ of bounded deterministic functions, a similar calculation of the quadratic variation of $\int_0^t (d\widehat{Z}(s) -$ $\sum_{j=k}^{N} g_j(s) \, dY_j(s))$, and the fact that the quadratic variation is an increasing process, shows that, for a.e. $s \in [0, t_1]$,

$$b_{ZZ}(s) - 2\sum_{j=k}^{N} g_j(s) b_{Zj}(s) + \sum_{j=k}^{N} g_j(s)^2 b_{jj}(s) \geq 0.$$

Taking $g_j(s) = \min(b_{Zj}(s)/b_{jj}(s), M)$, $j = 1, \ldots, N$, and then letting $M \to \infty$ and $N \to \infty$ yields $b_{ZZ}(s) \geq \sum_{j=k}^{\infty} b_{Zj}(s)^2/b_{jj}(s)$ for a.e. $s \in [0, t_1]$. Combined with (5.53), this shows that, for a.e. $s \in [0, t_1]$,

$$(5.54) \qquad 4\lambda e^{2s} = b_{ZZ}(s) = \sum_{j=k}^{\infty} b_{Zj}(s)^2 / b_{jj}(s) = 4 e^{3s} \sum_{j=k}^{\infty} q'_j(e^{-s}).$$

It is easily seen from (2.9) and Condition 2.1(ii) that $\sum_{j=k}^{\infty} q'_j(p)$ can be expressed as a convergent power series in $p \in [0, 1]$; thus, (5.54) holds for all $s \geq 0$ and yields

$$\sum_{j=k}^{\infty} q'_j(p) = \lambda p, \qquad 0 \leq p \leq 1.$$

By integration, since $q_j(0) = 0$ for $j > 0$ by (2.4), this gives

$$(5.55) \qquad \sum_{j=k}^{\infty} q_j(p) = \frac{\lambda}{2} p^2, \qquad 0 \leq p \leq 1.$$

As $p \to 0$, the left-hand side is $O(p^k)$, so (5.55) requires $k = 2$. Moreover, by (2.4) again,

$$\frac{\lambda}{2} p^2 = \sum_{j=2}^{\infty} \mathbb{P}(D_p \geq j) = \mathbb{E} D_p - \mathbb{P}(D_p \geq 1) = p\lambda - 1 + \mathbb{P}(D_p = 0),$$

and thus, for $q = 1 - p \in [0, 1]$,

$$\sum_{j=0}^{\infty} p_j q^j = \mathbb{P}(D_{1-q} = 0) = \frac{\lambda}{2}(1-q)^2 - \lambda(1-q) + 1 = 1 - \frac{\lambda}{2} + \frac{\lambda}{2} q^2.$$

Consequently, $p_j = 0$ for $j \neq 0, 2$, which contradicts our assumption $\mathbb{P}(D > k) > 0$. (Although somewhat hidden in the argument above, the conceptual reason behind the proof is that $W_n$ will jump without any change in $V_{rn}$, $r \geq k$, every time a light white ball dies, and this occurs quite often.) □

In the above covariance formulae, $dq_j(p) = q'_j(p) \, dp$, which can be expressed in $(p_l)_l$ using (2.9). In particular, (5.39) and (2.9) yield

$$\text{Cov}(\widehat{X}_r(-\ln x), \widehat{Z}(-\ln x))$$

$$= 2 \sum_{j=r}^{\infty} \sum_{l=j}^{\infty} l p_l \binom{l-1}{j-1} \binom{j-1}{r-1}$$

$$\times \int_x^1 (p-x)^{j-r} p^{-j-2} p^{j-1} (1-p)^{l-j}\, dp$$

$$= 2 \sum_{l=r}^{\infty} p_l \sum_{j=r}^{l} l \binom{l-1}{r-1} \binom{l-r}{j-r} \int_x^1 (p-x)^{j-r} (1-p)^{l-j} p^{-3}\, dp$$

$$= 2 \sum_{l=r}^{\infty} p_l r \binom{l}{r} \int_x^1 (1 - p + p - x)^{l-r} p^{-3}\, dp$$

$$= r \sum_{l=r}^{\infty} p_l \binom{l}{r} (1-x)^{l-r} (x^{-2} - 1).$$

Consequently, (5.40) simplifies to

(5.56)  $$\sigma_{rW}(x) = r x^r (1 - x^2) \sum_{l=r}^{\infty} p_l \binom{l}{r} (1-x)^{l-r}.$$

Unfortunately, we have not found any similar simplification for, for example, $\sigma_{rr}(x)$.

EXAMPLE 5.2. In the special case where $D \sim \mathrm{Po}(\lambda)$, we have by Lemma 2.4 and the line following (1.2)

$$q_j'(p) = \lambda \psi_j'(\lambda p) = \lambda \pi_{j-1}(\lambda p) = e^{-\lambda p} \frac{\lambda^j p^{j-1}}{(j-1)!}$$

and by (5.56) [or from (5.40)],

$$\sigma_{rW}(x) = r x^r (1 - x^2) \sum_{l=r}^{\infty} e^{-\lambda} \frac{\lambda^l}{(l-r)! r!} (1-x)^{l-r}$$

$$= x^r (1 - x^2) \frac{e^{-\lambda} \lambda^r}{(r-1)!} \sum_{i=0}^{\infty} \frac{(\lambda(1-x))^i}{i!} = \frac{(\lambda x)^r e^{-\lambda x}}{(r-1)!} (1 - x^2),$$

but we see no significant simplification of (5.41).

**6. Proof of Theorems 3.4 and 3.5 for $G^*(n, (d_i)_1^n)$.** We will now use Theorem 3.1 to prove Theorems 3.4 and 3.5 for the random multigraph $G^*(n, (d_i)_1^n)$. The simple random graph $G(n, (d_i)_1^n)$ will be treated in the next section. As a preparation, the reader might (on second reading) observe that the proofs below hold also if we allow $O(1)$ golden edges, as in Theorem 3.1, now interpreting $\mathrm{Core}_k$ as the remainder (excluding isolated vertices) when we stop, provided that in Theorem 3.5 we replace "Core$_k = \varnothing$" by "$H_n(\tau_n) < \bar{\delta}n$" for some fixed, sufficiently small, $\bar{\delta} > 0$.

PROOF OF THEOREM 3.4. Note that $l(0) = l_n(0) = 0$, so $\hat{p}$ and $\hat{p}_n$ are well-defined numbers in $[0,1]$. Moreover, if $\delta > 0$ is small enough, then $l(\hat{p} - \delta) < 0$ and $l(\hat{p} + \delta) > 0$. Since $l_n \to l$ by Lemma 2.3, this implies $l_n(\hat{p} - \delta) < 0$ and $l_n(\hat{p} + \delta) > 0$ for large enough $n$, and thus, $l_n$ then has a zero $\hat{p}_n$ in $(\hat{p} - \delta, \hat{p} + \delta)$. Furthermore, the uniform convergence of $l_n$ to $l$ given by Lemma 2.3, and the fact that $l > 0$ on the compact interval $[\hat{p} + \delta, 1]$, imply that if $n$ is large enough, then $l_n > 0$ on $[\hat{p} + \delta, 1]$. It follows that, for large enough $n$, $\hat{p}_n \in (\hat{p} - \delta, \hat{p} + \delta)$. Since $\delta$ can be chosen arbitrarily small, this shows $\hat{p}_n \to \hat{p}$. We change variables and define $\hat{t} := -\ln \hat{p}$ and $t_n := -\ln \hat{p}_n$; thus, $t_n \to \hat{t}$.

It was shown in [17], proof of Theorem 2.3, that $\tau_n \xrightarrow{\mathrm{P}} \hat{t}$. [This was proved using methods similar to those used here, as a simple consequence of (3.9), and extends to the case with golden balls.] Consequently, we can apply Theorem 3.1 with $\tau'_n = \tau_n$ and $t_0 = \hat{t}$. We simplify, at least conceptually, by using the Skorokhod coupling theorem ([19], Theorem 4.30) which shows that we can assume that all random variables are defined on the same probability space and that both $\tau_n \to \hat{t}$ and the limits (3.5)–(3.7) hold a.s. Since the limits $Z_B$, $Z_H$ and $Z_L$ are a.s. continuous, this means that a.s. (3.5)–(3.7) hold uniformly on the interval $[0, \hat{t} + 1]$. Taking $t = \tau_n$ (which a.s. is less that $\hat{t} + 1$ for large $n$), we see that, in particular, a.s.

$$\begin{aligned}(6.1) \quad L_n(\tau_n) &= n\check{l}_n(\tau_n) + n^{1/2} Z_L(\tau_n \wedge \hat{t}) + o(n^{1/2}) \\ &= n\check{l}_n(\tau_n) + n^{1/2} Z_L(\hat{t}) + o(n^{1/2}),\end{aligned}$$

by the continuity of $Z_L$. Thus, since $L_n(\tau_n) = -1$,

$$(6.2) \quad \check{l}_n(\tau_n) = -n^{-1/2} Z_L(\hat{t}) + o(n^{-1/2}).$$

Moreover, by the mean-value theorem, for some $\xi_n$ in the interval $[t_n, \tau_n]$ or $[\tau_n, t_n]$,

$$(6.3) \quad \check{l}_n(\tau_n) = \check{l}_n(\tau_n) - \check{l}_n(t_n) = \check{l}'_n(\xi_n)(\tau_n - t_n).$$

As $n \to \infty$, we have $\tau_n \to \hat{t}$ and $t_n \to \hat{t}$, and thus, $\xi_n \to \hat{t}$. Hence, the uniform convergence of $\check{l}'_n$ (see Lemma 2.3) implies $\check{l}'_n(\xi_n) \to \check{l}'(\hat{t}) = -e^{-\hat{t}} l'(\hat{p}) = -\hat{p}\alpha$. It follows by (6.3) and (6.2) that a.s.

$$\tau_n - t_n = \left(-\frac{1}{\hat{p}\alpha} + o(1)\right) \check{l}_n(\tau_n) = n^{-1/2} \frac{1}{\alpha \hat{p}} (Z_L(\hat{t}) + o(1)).$$

Consequently, using the mean-value theorem again and the analogue of (6.1) for $H_n$, a.s., for some $\xi'_n \to \hat{t}$,

$$
\begin{aligned}
n^{-1/2} H_n(\tau_n) &= n^{1/2} \check{h}_n(\tau_n) + Z_H(\hat{t}) + o(1) \\
&= n^{1/2} \check{h}_n(t_n) + n^{1/2} \check{h}'_n(\xi'_n)(\tau_n - t_n) + Z_H(\hat{t}) + o(1) \\
&= n^{1/2} h_n(\hat{p}_n) + \check{h}'(\hat{t}) \frac{1}{\alpha \hat{p}} Z_L(\hat{t}) + Z_H(\hat{t}) + o(1) \\
&= n^{1/2} h_n(\hat{p}_n) - \frac{h'(\hat{p})}{\alpha} Z_L(\hat{t}) + Z_H(\hat{t}) + o(1),
\end{aligned}
\tag{6.4}
$$

and similarly for $B_n$. The result (3.10) follows since $e(\text{Core}_k) = \frac{1}{2} H_n(\tau_n)$ and $v(\text{Core}_k) = B_n(\tau_n)$.

Finally, it is easily seen that if $\mathbb{P}(D > k) = 0$, then $h(p) = k p_k p^k$ and either $l(1) = 0$ and $\hat{p} = 1$, or $l(p) > 0$ for $0 < p \leq 1$ and $\hat{p} = 0$. But we have assumed that $0 < \hat{p} < 1$, so $\mathbb{P}(D > k) > 0$ and the nonsingularity of the covariance matrix follows from the nonsingularity in Theorem 3.1. □

PROOF OF THEOREM 3.5. Let $0 < \delta_0 < \hat{p}$, and let $\varepsilon > 0$ be so small that $\delta_0 < \hat{p} - \varepsilon < \hat{p} < \hat{p} + \varepsilon < 1$ and $l'' > \beta/2$ on $I_\varepsilon := [\hat{p} - \varepsilon, \hat{p} + \varepsilon]$. Since $l''_n \to l''$ uniformly, it follows that $l''_n > \beta/4 > 0$ on $I_\varepsilon$, provided $n$ is large. Hence, for such $n$, $l_n$ has a unique minimum point $\bar{p}_n$ in $I_\varepsilon$. Moreover, since $l_n \to l$ uniformly, $l_n(\bar{p}_n) = \min_{I_\varepsilon} l_n \to \min_{I_\varepsilon} l = 0$. On the other hand, $\eta := \min_{[\delta_0, 1] \setminus I_\varepsilon} l > 0$, so if $n$ is large enough, then $\min_{[\delta_0, 1] \setminus I_\varepsilon} l_n > \eta/2 > \min_{I_\varepsilon} l_n = l_n(\bar{p}_n)$. Consequently, for large $n$, $\bar{p}_n$ is the unique minimum point of $l_n$ in $[\delta_0, 1]$. In the sequel we consider only $n$ such that such a unique minimum point exists, and we redefine $\bar{p}_n$ to be this minimum point. We have shown that the two definitions give the same result for large $n$. In particular, $\bar{p}_n \in I_\varepsilon$ for large $n$. Since $\varepsilon$ can be chosen arbitrarily small, this shows that $\bar{p}_n \to \hat{p}$.

We change variables and define $\hat{t} := -\ln \hat{p}$ and $t_n := -\ln \bar{p}_n \to \hat{t}$. Thus, $\check{l}(\hat{t}) = 0$, but $\check{l}(t) > 0$ for $t \geq 0$ with $t \neq \hat{t}$. Suppose that $T_1$ and $T_2$ are any fixed numbers with $\hat{t} < T_1 < T_2$. Then $\eta_1 := \min\{\check{l}(t) : t \in [T_1, T_2]\} > 0$. By the uniform convergence of $\check{l}_n \to \check{l}$, $\min_{[T_1, T_2]} \check{l}_n > \eta_1/2$ for large $n$. Hence, if $T_1 \leq \tau_n \leq T_2$, and $n$ is large enough, then $\check{l}_n(\tau_n) - L_n(\tau_n)/n = \check{l}_n(\tau_n) + 1/n > \eta_1/2$, because $L_n(\tau_n) = -1$. On the other hand, it follows from (3.9) that $\mathbb{P}(\check{l}_n(\tau_n) - L_n(\tau_n)/n > \eta_1/2) \to 0$. Consequently,

$$\mathbb{P}(T_1 \leq \tau_n \leq T_2) \to 0. \tag{6.5}$$

Similarly, if $T < \hat{t}$, then

$$\mathbb{P}(0 \leq \tau_n \leq T) \to 0. \tag{6.6}$$

By [17], Lemma 5.1 (a version of a classical result by Łuczak [21]), there exists $\delta > 0$ such that whp $G^*(n, (d_i)_1^n)$ has no nonempty $k$-core with fewer that $\delta n$ vertices. Choose $T_2$ such that $\lambda e^{-2T_2} < \delta$, and let $\delta_1 := \delta - \lambda e^{-2T_2} > 0$.

Since $W_n(0)/n = (2m-1)/n \to \lambda$, it follows by (5.4) that $|\frac{1}{n}W_n(T_2 \wedge \tau_n) - \lambda e^{-2(T_2 \wedge \tau_n)}| \xrightarrow{\mathrm{P}} 0$, and thus, whp, if $\tau_n \geq T_2$, then $|\frac{1}{n}W_n(T_2) - \lambda e^{-2T_2}| < \delta_1$, and thus, $\frac{1}{n}W_n(T_2) < \lambda e^{-2T_2} + \delta_1 = \delta$. Hence, if $\tau_n \geq T_2$, then whp there are less than $\delta n$ white balls remaining at $T_2$, and thus even fewer at $\tau_n$, so the $k$-core has fewer that $\delta n$ vertices. By the result just quoted, this shows that if $\tau_n \geq T_2$, then whp the $k$-core is empty. In conjunction with (6.5), this implies that, for every fixed $T_1 > \hat{t}$,

(6.7) $\qquad \mathbb{P}(\mathrm{Core}_k \neq \varnothing \text{ and } \tau_n \geq T_1) \to 0.$

We next replace the fixed $T_1$ by a sequence $t'_n \to \hat{t}$. We claim we can define $t'_n \geq \max\{\hat{t}, t_n\}$ in such a way that $t'_n \to \hat{t}$ and

(6.8) $\qquad \mathbb{P}(\mathrm{Core}_k \neq \varnothing \text{ and } \tau_n \geq t'_n) \to 0.$

Let us do so explicitly. For every $i \in \mathbb{N}$, we may define $n_i = \max\{n_{i-1}+1, n'_i\}$, where $n'_i$ is the minimum natural number such that, for all $n \geq n'_i$,

$$\mathbb{P}(\mathrm{Core}_k \neq \varnothing \text{ and } \tau_n \geq \hat{t} + 1/i) \leq 1/i.$$

[It is clear by (6.7) that $n_i$ exists.] The numbers $n_i$ form an increasing sequence. Now define $t''_n = \hat{t} + 1/i$ for each $n_i \leq n < n_{i+1}$. Then for $n \geq n_i$,

$$\mathbb{P}(\mathrm{Core}_k \neq \varnothing \text{ and } \tau_n \geq t''_n) \leq \frac{1}{i}.$$

Further, $t''_n \to \hat{t}$. Finally, we can let $t'_n = \max\{t''_n, t_n\}$, and this is a sequence with required properties.

We now define $\tau'_n := \tau_n \wedge t'_n$. If $\tau'_n = t'_n$, then $\tau_n \geq t'_n$ and thus, by (6.8), whp the $k$-core is empty. Conversely, if $\tau'_n < t'_n$, then $\tau_n = \tau'_n < t'_n < \hat{t} + 1$ for large $n$. Let $\eta_2 := \lambda e^{-2(\hat{t}+2)}$. By (5.4) again, if $\tau_n < \hat{t} + 1$, then whp $n^{-1}W_n(\tau_n) > \lambda e^{-2\tau_n} - \eta_2 > 0$, and thus, there are some white balls left at $\tau_n$, which shows that there is a nonempty $k$-core. Consequently, whp the $k$-core is empty if and only if $\tau'_n = t'_n$.

Note further that $\tau'_n \leq t'_n$ by definition, and $t'_n \to \hat{t}$. Moreover, by (6.6), for every $T < \hat{t}$, whp $\tau_n \geq T$ and, thus, $\tau'_n \geq T$. Consequently, $\tau'_n \xrightarrow{\mathrm{P}} \hat{t}$.

We can thus apply Theorem 3.1, with $t_0 = \hat{t}$. Since $0 < \hat{p} < 1$, the non-singularity statement in Theorem 3.1 applies as in the proof of Theorem 3.4; in particular, $\sigma^2 = \sigma_{LL}^2 > 0$.

As in the proof of Theorem 3.4, we use the Skorohod coupling theorem [19], Theorem 4.30 and assume that the limits (3.5)–(3.7) and $\tau'_n \to \hat{t}$ hold

a.s. In particular, (3.7) then yields a.s., using $\tau'_n \to \hat{t}$ and the continuity of $Z_L$,

(6.9) $$n^{-1/2}L_n(\tau'_n) = n^{1/2}\breve{l}_n(\tau'_n) + Z_L(\hat{t}) + o(1).$$

If $\tau'_n < t'_n$, and thus $\tau_n = \tau'_n$ and $L_n(\tau'_n) = L_n(\tau_n) = -1$, then (6.9) yields

(6.10) $$n^{1/2}\breve{l}_n(\tau'_n) = -Z_L(\hat{t}) + o(1),$$

and, since $\bar{p}_n$ is the minimum point on an interval including $\tau_n$,

(6.11) $$n^{1/2}l_n(\bar{p}_n) = n^{1/2}\breve{l}_n(t_n) \leq n^{1/2}\breve{l}_n(\tau'_n) = -Z_L(\hat{t}) + o(1).$$

In particular, for any $\varepsilon > 0$, if $\text{Core}_k \neq \varnothing$, then whp $\tau'_n < t'_n$ and, by (6.11), whp

(6.12) $$n^{1/2}l_n(\bar{p}_n) \leq -Z_L(\hat{t}) + \varepsilon.$$

Conversely, if $\tau'_n = t'_n$, then a.s. $\tau_n > t'_n \geq t_n$ (since $\tau_n$ has a continuous distribution) so $L_n(t_n) \geq 0$. Moreover, (3.7) yields

(6.13) $$0 \leq n^{-1/2}L_n(t_n) = n^{1/2}\breve{l}_n(t_n) + Z_L(\hat{t}) + o(1).$$

Consequently,

(6.14) $$n^{1/2}l_n(\bar{p}_n) = n^{1/2}\breve{l}_n(t_n) \geq -Z_L(\hat{t}) + o(1).$$

In particular, for any $\varepsilon > 0$, if $\text{Core}_k = \varnothing$, then whp $\tau'_n = t'_n$ and, by (6.14), whp

(6.15) $$n^{1/2}l_n(\bar{p}_n) \geq -Z_L(\hat{t}) - \varepsilon.$$

In the case (i) of the theorem, (6.12) fails for large $n$ (and any $\varepsilon > 0$), and thus, whp $\text{Core}_k = \varnothing$. Similarly, in case (iii), (6.15) fails for large $n$ (and any $\varepsilon > 0$), and thus, whp $\text{Core}_k \neq \varnothing$.

In case (ii), it follows similarly that if $Z_L(\hat{t}) > -\zeta + \varepsilon$, then (6.12) fails for large $n$, and thus, whp $\text{Core}_k = \varnothing$, while if $Z_L(\hat{t}) < -\zeta - \varepsilon$, then (6.15) fails for large $n$ and thus whp $\text{Core}_k \neq \varnothing$. Since $\varepsilon > 0$ is arbitrary, this means that whp $\text{Core}_k \neq \varnothing \iff Z_L(\hat{t}) < -\zeta$, and thus, since $Z_L(\hat{t}) \sim N(0, \sigma^2)$ by (3.8),

$$\mathbb{P}(\text{Core}_k \neq \varnothing) \to \mathbb{P}(Z_L(\hat{t}) < -\zeta) = \Phi(-\zeta/\sigma),$$

as asserted.

Assume now that $\tau_n = \tau'_n < t'_n$ so that the $k$-core is nonempty, and let us localize $\tau_n$ more precisely. Note first that, since $l'(\hat{p}) = 0$ and $l''(\hat{p}) = \beta > 0$, it follows that $\breve{l}'(\hat{t}) = 0$ and $\breve{l}''(\hat{t}) = \breve{\beta} := \beta\hat{p}^2 > 0$.

Using Taylor's formula and the fact that $\breve{l}'_n(t_n) = 0$, for some $\xi_n$ between $\tau'_n$ and $t_n$,

$$\breve{l}_n(\tau'_n) = \breve{l}_n(t_n) + \tfrac{1}{2}\breve{l}''_n(\xi_n) \cdot (\tau'_n - t_n)^2,$$

and thus by (6.10), recalling $\check{l}_n(t_n) = l_n(\bar{p}_n)$,

(6.16) $\quad \frac{1}{2}\check{l}''_n(\xi_n) \cdot (\tau'_n - t_n)^2 = n^{-1/2}(-Z_L(\hat{t}) - n^{1/2}l_n(\bar{p}_n) + o(1)).$

Since $\tau'_n \to \hat{t}$ and $t_n \to \hat{t}$, we have $\xi_n \to \hat{t}$ a.s., and thus, $\check{l}''_n(\xi_n) \to \check{l}''(\hat{t}) = \check{\beta} > 0$ a.s.

In case (ii) we thus have, when $Z_L(\hat{t}) < -\zeta$,

(6.17) $\quad\quad\quad (\tau'_n - t_n)^2 = 2\check{\beta}^{-1}n^{-1/2}(-Z_L(\hat{t}) - \zeta + o(1)).$

If $\tau'_n > t_n$, then $L_n(t_n) \geq 0$ and (6.13) holds, which by (6.16) implies

(6.18) $\quad\quad\quad \frac{1}{2}\check{l}''_n(\xi_n) \cdot (\tau'_n - t_n)^2 \leq o(n^{-1/2}),$

which contradicts (6.17) for large $n$. Consequently, when $Z_L(\hat{t}) < -\zeta$, (6.17) yields

(6.19) $\quad \tau_n = \tau'_n = t_n - n^{-1/4}((2/\check{\beta})^{1/2}(-Z_L(\hat{t}) - \zeta)^{1/2} + o(1)).$

We now obtain from (3.6) [cf. (6.4)], recalling $\check{h}'_n(t_n) \to \check{h}'(\hat{t}) = -\hat{p}h'(\hat{p})$,

$$\begin{aligned}H_n(\tau_n) &= n\check{h}_n(\tau_n) + n^{1/2}(Z_H(\hat{t}) + o(1)) \\ (6.20) \quad &= n\check{h}_n(t_n) + n(\check{h}'_n(t_n) + o(1))(\tau_n - t_n) + O(n^{1/2}) \\ &= nh_n(\bar{p}_n) + n^{3/4}(h'(\hat{p})(2/\beta)^{1/2}(-Z_L(\hat{t}) - \zeta)^{1/2} + o(1))\end{aligned}$$

and similarly for $B_n$. The result follows, with $Z := -Z_L(\hat{t})/\sigma \sim N(0,1)$, since $e(\mathrm{Core}_k) = \frac{1}{2}H_n(\tau_n)$ and $v(\mathrm{Core}_k) = B_n(\tau_n)$. (Note that $Z_H$ and $Z_B$ give insignificant contributions in this case.)

In case (iii), we similarly obtain from (6.16)

$$(\tau'_n - t_n)^2 = 2\check{\beta}^{-1}|l_n(\bar{p}_n)|(1 + o(1))$$

and, since by (6.18) we can again exclude $\tau'_n > t_n$ for large $n$,

(6.21) $\quad \tau_n - t_n = \tau'_n - t_n = -(2/\check{\beta})^{1/2}|l_n(\bar{p}_n)|^{1/2}(1 + o(1)).$

The random fluctuations here turn out to be of a smaller order. To capture them, we note that, for large $n$, using the argument at the start of the proof [since now $l_n(\bar{p}_n) < 0$], that $l_n(p) = 0$ for exactly two values of $p$ in $I_\varepsilon$, say, $\hat{p}_n$ and $\hat{p}'_n$, with $\hat{p} - \varepsilon < \hat{p}'_n < \bar{p}_n < \hat{p}_n < \hat{p} + \varepsilon$, and not for any other $p \in [\delta_0, 1]$. Let $\hat{t}_n := -\ln\hat{p}_n < t_n$, so $\check{l}_n(\hat{t}_n) = 0$. By Taylor's formula as for (6.16), also

(6.22) $\quad\quad\quad \hat{t}_n - t_n = -(2/\check{\beta})^{1/2}|l_n(\bar{p}_n)|^{1/2}(1 + o(1)).$

Finally, by Taylor's formula again, for some $\xi_n$ between $\tau_n$ and $\hat{t}_n$,

(6.23) $\quad\quad\quad \check{l}_n(\tau_n) = \check{l}_n(\tau_n) - \check{l}_n(\hat{t}_n) = \check{l}'_n(\xi_n)(\tau_n - \hat{t}_n).$

Since $\check{l}'_n(t_n) = 0$, we similarly have for some $\xi'_n$ between $\xi_n$ and $t_n$,

(6.24) $$\check{l}'_n(\xi_n) = \check{l}'_n(\xi_n) - \check{l}'_n(t_n) = (\xi_n - t_n)\check{l}''_n(\xi'_n).$$

From (6.21) and (6.22), $\tau_n - t_n \sim \hat{t}_n - t_n$, and so also $\xi_n - t_n \sim \hat{t}_n - t_n$. Combining this with (6.24), (6.22), as well as the fact that $\check{l}''_n(\xi'_n) \to \check{l}''(\hat{t}) = \check{\beta}$, we obtain

$$\check{l}'_n(\xi_n) = \check{\beta}(\hat{t}_n - t_n)(1 + o(1)) = -(2\check{\beta})^{1/2}|l_n(\bar{p}_n)|^{1/2}(1 + o(1)).$$

Hence, (6.23) and (6.10) yield

$$\tau_n - \hat{t}_n = |l_n(\bar{p}_n)|^{-1/2} n^{-1/2} (2\check{\beta})^{-1/2}(Z_L(\hat{t}) + o(1)).$$

The result now follows by (3.5) and (3.6), analogously to (6.20), with $Z' := -(2\beta)^{-1/2} Z_L(\hat{t})$.

For the final statements, from Taylor's formula and $l'_n(\bar{p}_n) = 0$,

(6.25) $$l'_n(\hat{p}) = (\beta + o(1))(\hat{p} - \bar{p}_n),$$

and hence

$$l_n(\hat{p}) - l_n(\bar{p}_n) = \tfrac{1}{2}(\beta + o(1))(\hat{p} - \bar{p}_n)^2 = ((2\beta)^{-1} + o(1))l'_n(\hat{p})^2,$$

so, if $l'_n(\hat{p})^2 = o(l_n(\hat{p}))$, then $l_n(\bar{p}_n)/l_n(\hat{p}) \to 1$, which yields the first claim. The second follows from (6.25) and $b_n(\bar{p}_n) - b_n(\hat{p}) = O(|\hat{p} - \bar{p}_n|)$, and similarly for $h_n$. □

**7. The simple graph $G(n, (d_i)_1^n)$.** We have proved Theorems 3.4 and 3.5 for the random multigraph $G^*(n, (d_i)_1^n)$, and the next step is to show that they hold for the random simple graph $G(n, (d_i)_1^n)$ too, that is, that they hold for $G^*(n, (d_i)_1^n)$ conditioned on this random multigraph being simple. For results that can be stated as convergence in probability, or stating that some event holds whp, this transfer is immediate: it suffices to have that

(7.1) $$\liminf \mathbb{P}(G^*(n, (d_i)_1^n) \text{ is simple}) > 0.$$

(As is well known (see, e.g., [22]), (7.1) holds under Condition 2.1(ii); see [16] for a general necessary and sufficient condition.) For example, this holds in our previous paper [17]. However, with distributional results, such as those in the present paper, the transfer to $G(n, (d_i)_1^n)$ is much more delicate, and amounts to showing that the variables we study are asymptotically independent of the event that $G^*(n, (d_i)_1^n)$ is simple.

Our basic tool will be the following general probabilistic result. Recall from Section 2 that an indicator (of some event) is a random variable taking values in the set $\{0, 1\}$.

PROPOSITION 7.1. *Assume that $X_n$, $n \geq 1$, and $X$ are random variables with values in some metric space $\mathcal{S}$ and that, for each $n \geq 1$, $(I_{n\alpha})_{\alpha \in \mathcal{A}_n}$, is a (finite) family of indicator random variables defined on the same probability space as $X_n$. Let $W_n := \sum_{\alpha \in \mathcal{A}_n} I_{n\alpha}$ and let $\mathcal{E}_n$ be the event*

$$\mathcal{E}_n := \{I_{n\alpha} = 0 \text{ for every } \alpha \in \mathcal{A}_n\} = \{W_n = 0\}.$$

*Assume further the following, as $n \to \infty$:*

(i) $X_n \xrightarrow{d} X$.

(ii) *For any fixed $\ell \geq 1$, if we for each $n$ select $\ell$ indices $\alpha_1^n, \ldots, \alpha_\ell^n \in \mathcal{A}_n$ such that $\mathbb{P}(I_{n\alpha_1^n} = \cdots = I_{n\alpha_\ell^n} = 1) > 0$, then*

(7.2) $$(X_n \mid I_{n\alpha_1^n} = \cdots = I_{n\alpha_\ell^n} = 1) \xrightarrow{d} X.$$

(iii) $W_n \xrightarrow{d} W$, *where $W$ is a random variable such that $\mathbb{P}(W = 0) > 0$ and the distribution of $W$ is determined by its moments.*

(iv) $\limsup_{n \to \infty} \mathbb{E}(W_n^j) < \infty$ *for every $j \geq 1$.*

*Then $(X_n \mid \mathcal{E}_n) \xrightarrow{d} X$.*

REMARK 7.2. The careful reader may observe that in (ii) it is conceivable that, for some $n$, it is impossible to select $\alpha_1^n, \ldots, \alpha_\ell^n \in \mathcal{A}_n$ satisfying the condition. However, we may just ignore such $n$. Indeed, since we may take $\alpha_1^n = \cdots = \alpha_\ell^n$, this happens if and only if $W_n = 0$ a.s., and if this happens for more than a finite number of values of $n$, then necessarily $W = 0$ a.s. by (iii), so $\mathbb{P}(\mathcal{E}_n) \to 1$ and the conclusion becomes trivial.

PROOF OF PROPOSITION 7.1. Recall that $X_n \xrightarrow{d} X$ is equivalent to $\mathbb{P}(X_n \in A) \to \mathbb{P}(X \in A)$ for all measurable sets $A \subseteq \mathcal{S}$ such that $\mathbb{P}(X \in \partial A) = 0$; see [2, Theorem 2.1].

Fix such a set $A$; we thus want to prove that $\mathbb{P}(X_n \in A \mid \mathcal{E}_n) \to \mathbb{P}(X \in A)$. By (7.2), for each $\ell \geq 1$,

(7.3) $$\mathbb{P}(X_n \in A \mid I_{n\alpha_1^n} = \cdots = I_{n\alpha_\ell^n} = 1) \xrightarrow{d} \mathbb{P}(X \in A)$$

for every choice of $\alpha_1^n, \ldots, \alpha_\ell^n \in \mathcal{A}_n$ such that $\mathbb{P}(I_{n\alpha_1^n} = \cdots = I_{n\alpha_\ell^n} = 1) > 0$. Moreover, (7.3) holds uniformly in all such choices, since otherwise we could for each $n$ select $\alpha_1^n, \ldots, \alpha_\ell^n$ such that the difference between the two sides in (7.3) is maximal (for this $n$) and obtain a contradiction.

Assume first that $\mathbb{P}(X \in A) > 0$. Then $\mathbb{P}(X_n \in A) > 0$ for large $n$; we consider such $n$ only. Then, by (7.3) and assumption (i),

$$\mathbb{E}(I_{n\alpha_1} \cdots I_{n\alpha_\ell} \mid X_n \in A) = \frac{\mathbb{P}(X_n \in A \text{ and } I_{n\alpha_1} = \cdots = I_{n\alpha_\ell} = 1)}{\mathbb{P}(X_n \in A)}$$

$$= \frac{\mathbb{P}(X_n \in A \mid I_{n\alpha_1} = \cdots = I_{n\alpha_\ell} = 1)}{\mathbb{P}(X_n \in A)} \mathbb{P}(I_{n\alpha_1} = \cdots = I_{n\alpha_\ell} = 1)$$

$$= \frac{\mathbb{P}(X \in A) + o(1)}{\mathbb{P}(X \in A) + o(1)} \mathbb{P}(I_{n\alpha_1} = \cdots = I_{n\alpha_\ell} = 1)$$

$$= (1 + o(1)) \mathbb{E}(I_{n\alpha_1} \cdots I_{n\alpha_\ell}),$$

where $o(1) \to 0$ as $n \to \infty$, uniformly in all choices of $\alpha_1, \ldots, \alpha_\ell \in \mathcal{A}_n$ such that $\mathbb{E}(I_{n\alpha_1} \cdots I_{n\alpha_\ell}) > 0$. Furthermore, this holds also, trivially, if $\mathbb{E}(I_{n\alpha_1} \cdots I_{n\alpha_\ell}) = 0$. Consequently, we can sum over all $\alpha_1, \ldots, \alpha_\ell \in \mathcal{A}_n$ and obtain

(7.4) $$\mathbb{E}(W_n^\ell \mid X_n \in A) = (1 + o(1)) \mathbb{E}(W_n^\ell).$$

Now, assumptions (iii) and (iv) imply that $\mathbb{E}(W_n^\ell) \to \mathbb{E}(W^\ell)$ for every $\ell$, and thus, (7.4) yields $\mathbb{E}(W_n^\ell \mid X_n \in A) \to \mathbb{E}(W^\ell)$. By the method of moments, and (iii), this implies that $(W_n \mid X_n \in A) \xrightarrow{d} W$, and thus,

$$\mathbb{P}(W_n = 0 \mid X_n \in A) \to \mathbb{P}(W = 0).$$

As a consequence, using (i) and (iii) again,

$$\mathbb{P}(X_n \in A \mid \mathcal{E}_n) = \frac{\mathbb{P}(X_n \in A \text{ and } W_n = 0)}{\mathbb{P}(W_n = 0)}$$

$$= \frac{\mathbb{P}(W_n = 0 \mid X_n \in A) \mathbb{P}(X_n \in A)}{\mathbb{P}(W_n = 0)}$$

$$\to \frac{\mathbb{P}(W = 0) \mathbb{P}(X \in A)}{\mathbb{P}(W = 0)} = \mathbb{P}(X \in A).$$

This is precisely the result sought in the case $\mathbb{P}(X \in A) > 0$.

When $\mathbb{P}(X \in A) = 0$, then, trivially, by (i) and (iii),

$$\mathbb{P}(X_n \in A \mid \mathcal{E}_n) \le \frac{\mathbb{P}(X_n \in A)}{\mathbb{P}(\mathcal{E}_n)} \to \frac{\mathbb{P}(X \in A)}{\mathbb{P}(W = 0)} = 0 = \mathbb{P}(X \in A).$$

Consequently, $\mathbb{P}(X_n \in A \mid \mathcal{E}_n) \to \mathbb{P}(X \in A)$ for every measurable $A$ with $\mathbb{P}(X \in \partial A) = 0$, and thus, $(X_n \mid \mathcal{E}) \xrightarrow{d} X$. $\square$

REMARK 7.3. The conditions in Proposition 7.1 may be weakened in several different ways. First, we may allow some exceptional $\ell$-tuples of indices in (ii). More precisely, it suffices that there exists, for every $\ell$ and $n$, a set $\mathcal{B}_{n\ell} \subset \mathcal{A}_n^\ell$ such that (7.2) holds when $(\alpha_1^n, \ldots, \alpha_\ell^n) \notin \mathcal{B}_{n\ell}$, and further, for every fixed $\ell$, $\sum_{(\alpha_1^n, \ldots, \alpha_\ell^n) \in \mathcal{B}_{n\ell}} \mathbb{E}(I_{n\alpha_1^n} \cdots I_{n\alpha_\ell^n}) \to 0$ as $n \to \infty$. Second, we may replace (iii) and (iv) by the assumptions that $\liminf_{n \to \infty} \mathbb{P}(\mathcal{E}_n) > 0$ and that there exists $B < \infty$ such that $M_j := \limsup_{n \to \infty} \mathbb{E}(W_n^j) \le B^j j!$ for all $j \ge 1$.

(The latter could be weakened to the Carleman condition $\sum_j M_j^{-1/2j} = \infty$, see [10, Section 4.10].) Indeed, $M_1 < \infty$ implies that the sequence $(W_n)_n$ is tight, so every subsequence has a subsubsequence converging in distribution to some $W$ (that may depend on the subsubsequence), and these assumptions imply that this $W$ has to satisfy the conditions in (iii). Hence, $(X_n \mid \mathcal{E}_n) \xrightarrow{d} X$ holds along all such subsubsequences, which implies that it holds for the full sequence.

PROOF OF THEOREMS 3.4 AND 3.5 FOR $G(n, (d_i)_1^n)$. To prove Theorems 3.4 and 3.5 for $G(n, (d_i)_1^n)$, with the help of Proposition 7.1, we first check that the conclusions can be put in the form $X_n \xrightarrow{d} X$, where $X_n$ and $X$ take values in some metric space $\mathcal{S}$.

For Theorem 3.4 and Theorem 3.5(iii), this is immediate, with $\mathcal{S} = \mathbb{R}^2$.

For Theorem 3.5(ii), we first observe that (3.11) is equivalent to $\mathbf{1}[\text{Core}_k \neq \varnothing] \xrightarrow{d} \text{Be}(\Phi(-\zeta/\sigma))$, where $\text{Be}(q)$ denotes a Bernoulli distribution with parameter $q$; here we can take, for example, $\mathcal{S} = \mathbb{R}$. If this holds, we then take $\mathcal{S} = [-\infty, \infty)^2$ and note that the second assertion in Theorem 3.5(ii) can be written as

$$n^{-3/4}(v(\text{Core}_k) - b_n(\bar{p}_n)n, e(\text{Core}_k) - \tfrac{1}{2}h_n(\bar{p}_n)n) \xrightarrow{d} X,$$

as random variables in $[-\infty, \infty)^2$, for a certain random variable $X$ with a point mass $\Phi(\zeta/\sigma)$ at $(-\infty, -\infty)$.

Finally, Theorem 3.5(i) can also be written as a convergence in distribution of indicator variables (this time with a degenerate limit), but in this case, the transfer to $G(n, (d_i)_1^n)$ is actually immediate by the comments at the beginning of this section.

We make one more modification for Theorem 3.5: we choose a sufficiently small $\delta > 0$ and replace the condition $\text{Core}_k \neq \varnothing$ by $e(\text{Core}_k) > \delta n$; by [17], Lemma 5.1, and (7.1), these conditions are whp equivalent for both $G^*(n, (d_i)_1^n)$ and $G(n, (d_i)_1^n)$, so the modified theorem is equivalent to the original one.

In all cases we are thus in the setting of Proposition 7.1, with $X_n$ some functional of $G^*(n, (d_i)_1^n)$, and condition (i) satisfied, since it is just the statement that the statement in question holds for $G^*(n, (d_i)_1^n)$.

We let $\mathcal{A}_n := \mathcal{A}_{n1} \cup \mathcal{A}_{n2}$, where $\mathcal{A}_{n1}$ is the family of all (unordered) pairs of two half-edges at the same vertex, and $\mathcal{A}_{n2}$ is the family of all pairs $\{\{x_1, y_1\}, \{x_2, y_2\}\}$ of two disjoint pairs of half-edges with $x_1$ and $x_2$ belonging to one vertex and $y_1$ and $y_2$ belonging to another. For $\alpha \in \mathcal{A}_n$, we let $I_{n\alpha}$ be the indicator of the event that the pair(s) in $\alpha$ occur in the random configuration. Then $W_n$ is the number of loops and pairs of parallel edges (other than loops) in $G^*(n, (d_i)_1^n)$, and $\mathcal{E}_n$ is the event that $G^*(n, (d_i)_1^n)$ is simple.

As is well known, Condition 2.1 implies that $W_n \xrightarrow{d} \mathrm{Po}(\widehat{\Lambda})$, with convergence of all moments, for some $\widehat{\Lambda} < \infty$; see, for example, [16], Theorem 7.1. [More precisely, $\widehat{\Lambda} = \Lambda + \Lambda^2$, where $\Lambda = \mathbb{E}\binom{D}{2}/\lambda$.] Consequently, conditions (iii) and (iv) in Proposition 7.1 are satisfied.

It remains to verify (ii). In our case, (ii) means that the results in Theorems 3.4 and 3.5 (modified as above) hold also if, for any fixed $\ell$, for every $n$, we select a set of at most $2\ell$ disjoint pairs of half-edges and condition the configuration on containing these pairs of half-edges. (Actually, we need this only for certain sets of pairs, but it is just as easy to prove it for arbitrary sets. Note that the modification in Theorem 3.5 is essential here, since, e.g., a set of $k$ parallel loops always yields a nonempty $k$-core.)

To see this, regard the half-edges in the selected pairs, as well as the edges in $G^*(n,(d_i)_1^n)$ corresponding to these pairs, as *special*. We consider two approximations of the $k$-core $\mathrm{Core}_k$ of $G^*(n,(d_i)_1^n)$. First, let $G^*(n,(d_i)_1^n)^-$ be the multigraph $G^*(n,(d_i)_1^n)$ with all special edges deleted, and let $\mathrm{Core}_k^-$ be the $k$-core of $G^*(n,(d_i)_1^n)^-$. Since $G^*(n,(d_i)_1^n)^-$ is a subgraph of $G^*(n,(d_i)_1^n)$, we have $\mathrm{Core}_k^- \subseteq \mathrm{Core}_k$. Second, apply the core-finding algorithm used in Section 3, but with all special balls colored golden and thus immune to death and recoloring, and let $\mathrm{Core}_k^+$ be the subgraph of $G^*(n,(d_i)_1^n)$ left at the end; thus, $\mathrm{Core}_k^+ \supseteq \mathrm{Core}_k$.

Since our variables $X_n$ are increasing functions of the numbers of vertices and edges in the $k$-core, it is sufficient to show that the results of Theorems 3.4 and 3.5 (modified) hold also for $\mathrm{Core}_k^-$ and $\mathrm{Core}_k^+$. For $\mathrm{Core}_k^+$, this follows, as remarked at the beginning of Section 6, by the proofs for $G^*(n,(d_i)_1^n)$ above, since the possibility of golden balls was included in Theorem 3.1. The results for $\mathrm{Core}_k^-$ follow from Theorems 3.4 and 3.5, since the random multigraph $G^*(n,(d_i)_1^n)^-$ is a random multigraph of the same type as $G^*(n,(d_i)_1^n)$, but with $O(1)$ vertex degrees $d_i$ decreased (with a total change of at most $4k$). Evidently, Condition 2.1 holds for the modified sequence too, so Theorems 3.4 and 3.5 apply.

This shows that condition (ii) in Proposition 7.1 holds in all cases, and hence its statement applies, showing that Theorems 3.4 and 3.5 hold for $G(n,(d_i)_1^n)$ too. $\square$

REMARK 7.4. Unfortunately, we have not been able to use Proposition 7.1 in a similar way to show that Theorem 3.1 holds for $G(n,(d_i)_1^n)$ too. The reason is that conditioning on the configuration containing even a single given pair of half-edges will, as far as we can see, mess up the process of balls and bins and introduce unwanted dependencies. And although we allow golden balls in Theorem 3.1, there is no monotonicity like the one that allowed us to consider $\mathrm{Core}_k^+$ and $\mathrm{Core}_k^-$ above. We discuss another attempt to extend Theorem 3.1 to $G(n,p)$ in Appendix; see also Remark 8.7.

**8. The random graphs $G(n,p)$ and $G(n,m)$.** We derive the results for $G(n,p)$ and $G(n,m)$ from our results for $G(n,(d_i)_1^n)$ by conditioning on the degree sequence. To see this in detail, let us be more general and consider a random (simple) graph $G_n$ with $n$ vertices labeled $1,\ldots,n$ and some random distribution of the edges such that any two graphs on $1,\ldots,n$ with the same degree sequence have the same probability of being attained by $G_n$. [Note that $G(n,p)$ and $G(n,m)$ are of this type.] Equivalently, conditioned on the degree sequence, $G_n$ is a random graph with that degree sequence of the type $G(n,(d_i)_1^n)$ introduced in Section 2. We may thus construct $G_n$ by first picking a random sequence $(d_i)_1^n = (d_i^{(n)})_1^n$ with the right distribution, and then choosing a random graph $G(n,(d_i)_1^n)$ for this $(d_i)_1^n$. (We assume that this is possible, which implies, in particular, that $\sum_{i=1}^n d_i$ is even.)

We will assume that Condition 2.1 holds in probability in the following sense. As usual, we let $u_r = u_r(n)$ be the (now random) number of vertices with degree $r$.

CONDITION 8.1. *Let $(d_i)_1^n = (d_i^{(n)})_1^n$ be the random sequence of vertex degrees of $G_n$. Then, for some probability distribution $(p_r)_{r=0}^\infty$ independent of $n$, with $p_0 < 1$:*

(i) *$u_r/n := \#\{i : d_i = r\}/n \xrightarrow{\mathrm{P}} p_r$ for every $r \geq 0$ as $n \to \infty$;*
(ii) *for every $A > 1$, we have $\sum_r u_r A^r = \sum_{i=1}^n A^{d_i} = O_{\mathrm{p}}(n)$.*

We begin with a technical lemma.

LEMMA 8.2. *If Condition 8.1 holds, we may, by replacing the random graphs $G_n$ by other random graphs $G'_n$ with the same distribution, assume that the random graphs are defined on a common probability space and that Condition 2.1 holds a.s.*

PROOF. From the Skorokhod coupling theorem ([19], Theorem 4.30) applied to the random sequences $(u_r(n))_{r=0}^\infty$, we may assume that the limit $u_r/n \to p_r$ in (i) holds a.s., for every $r \geq 0$. To "derandomize" (ii) as well, we refine that argument as follows. By (ii), for every $j \geq 1$ and $k \geq 1$, we may choose $C_{k,j}$ such that $\mathbb{P}(\sum_i k^{d_i} > C_{k,j} n) < 2^{-k}/j$. We may assume that, for every $k$, $C_{k,j}$ increases with $j$. Now condition on the event $\mathcal{E}_j := \{\sum_i k^{d_i} \leq C_{k,j} n$ for every $k \geq 1\}$, and note that $\mathbb{P}(\mathcal{E}_j) > 1 - 1/j$. Conditioned on $\mathcal{E}_j$, Condition 2.1(ii) holds uniformly. Let $\mathcal{E}_0 := \varnothing$. We apply the Skorokhod coupling theorem to $(u_r)_r$ conditioned on $\mathcal{E}_j \setminus \mathcal{E}_{j-1}$ for every $j \geq 1$ such that $\mathbb{P}(\mathcal{E}_j \setminus \mathcal{E}_{j-1}) > 0$; this shows that we can assume $u_r(n)/n \to p_r$ a.s. for every $r$ on $\mathcal{E}_j \setminus \mathcal{E}_{j-1}$, and we only have to combine these pieces for $j \geq 1$. □

We define the functions $b_n$, $h_n$ and $l_n$ as before, but now conditioned on the degree sequence $(d_i)_1^n$. Thus, we use (2.5)–(2.7) with $\mathbb{P}(X = l)$ and $\mathbb{E} X$ replaced by the random numbers $u_l/n$ and $2m/n = \sum_l l u_l/n$, respectively. Note that these functions are now random functions of $p$. [They depend on both $p$ and $(d_i)_1^n$.]

Define further the (deterministic) functions $b(p)$, $h(p)$, $l(p)$ as before, using (2.4)–(2.7) with $X = D$ having distribution $(p_r)_0^\infty$.

By Lemmas 8.2 and 2.3, $b_n \xrightarrow{\mathrm{P}} b$, $h_n \xrightarrow{\mathrm{P}} h$ and $l_n \xrightarrow{\mathrm{P}} l$, with convergence in probability, uniformly on $[0,1]$ and together with all derivatives.

THEOREM 8.3. *Let $G_n$ be as above and assume Condition 8.1.*

*If $k \geq 2$, then Theorem 3.4 holds also for the $k$-core of $G_n$; now $b_n$, $h_n$ and $l_n$ are random functions independent of $Z_B, Z_H, Z_L$; further, $\hat{p}_n$ is random, and $\hat{p}_n \xrightarrow{\mathrm{P}} \hat{p}$.*

*Similarly, if $k \geq 3$, then Theorem 3.5 holds for $G_n$ too, with the following modifications:* (a) *a unique minimum point $\bar{p}_n$ exists whp;* (b) *$\bar{p}_n$ is random, and $\bar{p}_n \xrightarrow{\mathrm{P}} \hat{p}$, $l_n(\bar{p}_n) \xrightarrow{\mathrm{P}} 0$;* (c) *$\hat{p}_n$ in* (iii) *exists whp (we can, if we like, use any supplementary definitions to have the random variables always defined);* (d) *in* (i), (ii) *and* (iii), *it suffices to have $n^{1/2} l_n(\bar{p}_n)$ convergent in probability, and $\zeta$ in* (ii) *may be random;* (e) *$Z$ should be independent of $\zeta$ and (3.11) has to be rewritten as $\mathbb{P}(\mathrm{Core}_k \neq \varnothing) \to \mathbb{P}(\sigma Z > \zeta)$.*

PROOF. The first statement (for $k \geq 2$) follows directly from Lemma 8.2 and Theorem 3.4. Note that $\hat{p}$, $\alpha$, $\hat{t}$, and the distribution of $(Z_B, Z_H, Z_L)$ do not depend on the random $(d_i)_1^n$, but only on the limit distribution $(p_r)_0^\infty$.

The second statement (for $k \geq 3$) follows similarly from Lemma 8.2 and Theorem 3.5, noting that if $n^{1/2} l_n(\bar{p}_n)$ converges in probability, we may include it in the application of the Skorohod coupling theorem in the proof of Lemma 8.2, and thus assume that it too converges a.s. For Theorem 3.5(ii), we also write the final conclusion as $(X_n \mid \mathrm{Core}_k \neq \varnothing) \xrightarrow{\mathrm{d}} (\mathbf{v}\sqrt{\sigma Z - \zeta} \mid \sigma Z > \zeta)$ for certain random vectors $X_n$ and a vector $\mathbf{v}$, and note that this and (3.11) together are equivalent to

$$\mathbb{E}(F(X_n); \mathrm{Core}_k \neq \varnothing) \to \mathbb{E}(F(\mathbf{v}\sqrt{\sigma Z - \zeta}); \sigma Z > \zeta)$$

for every bounded continuous function $F$ on $\mathbb{R}^2$. $\square$

We now specialize to $G(n,p)$ and $G(n,m)$, with $p \sim \lambda/n$ and $m \sim \lambda n/2$ for some fixed $\lambda > 0$, and observe first that Condition 8.1 holds and thus Theorem 8.3 applies.

LEMMA 8.4. *Consider either $G(n, \lambda_n/n)$ or $G(n, m)$, where $m = \lambda_n n/2$, and assume that $\lambda_n \to \lambda > 0$. Then, Condition 8.1 holds, with $p_r = \pi_r(\lambda)$ and, thus, $D \sim \mathrm{Po}(\lambda)$.*

PROOF. Condition 8.1(i) follows from (8.7) below, or by elementary estimates of mean and variance which we omit; cf. (8.10).

For (ii), similar elementary estimates show that if $M := \sup \lambda_n$, then
$$n^{-1} \mathbb{E}\, u_r(n) = O(M^r/r!)$$
uniformly in $r \geq 0$, and thus, for each $A$, $n^{-1} \mathbb{E} \sum_r u_r(n) A^r = O(1)$. □

In the cases of $G(n,p)$ and $G(n,m)$, the random functions $b_n$, $h_n$, $l_n$ are themselves asymptotically Gaussian processes. We give a precise statement.

THEOREM 8.5. *Let $b_n$, $h_n$, $l_n$ be the random functions defined above either for the random graph $G(n, \lambda_n/n)$ or for $G(n,m)$ with $m = \lambda_n n/2$, where $\lambda_n \to \lambda > 0$. Then, jointly in $D[0,1]$,*

(8.1) $$n^{1/2}(b_n(p) - b_{\text{Po}(\lambda_n)}(p)) \xrightarrow{d} U_B(p),$$

(8.2) $$n^{1/2}(h_n(p) - h_{\text{Po}(\lambda_n)}(p)) \xrightarrow{d} U_H(p),$$

(8.3) $$n^{1/2}(l_n(p) - l_{\text{Po}(\lambda_n)}(p)) \xrightarrow{d} U_L(p),$$

*where $U_B$, $U_H$ and $U_L$ are continuous Gaussian processes on $[0,1]$ with mean $0$ and covariances that satisfy, for $0 \leq p \leq 1$ and $\nu, \varkappa \in \{B, H, L\}$,*

(8.4) $$\text{Cov}(U_\nu(p), U_\varkappa(p)) = \sigma^*_{\nu\varkappa}(p),$$

*where $\sigma^*_{\nu\varkappa}$ are given by (8.13)–(8.23) below and either (8.8) for $G(n, \lambda_n/n)$ or (8.12) for $G(n,m)$. For $G(n,m)$, $U_L(p) = -U_H(p)$.*

PROOF. Consider first $G(n, \lambda_n/n)$. It is shown in [1] (see also [18], Example 6.35) that each $u_r = u_r(n)$ is asymptotically normal. More precisely,

(8.5) $$n^{-1/2}(u_r(n) - \mathbb{E}\, u_r(n)) \xrightarrow{d} U_r \sim N(0, \varphi_{rr}),$$

where, recalling the notation (1.1),
$$\varphi_{rr} := \lim_{n \to \infty} n^{-1} \text{Var}\, u_r(n) = \pi_r(\lambda)^2 \left(\frac{(r-\lambda)^2}{\lambda} - 1\right) + \pi_r(\lambda).$$

Moreover, with $p_n := \lambda_n/n$,

(8.6) $$\mathbb{E}\, u_r(n) = n \binom{n-1}{r} p_n^r (1-p_n)^{n-1-r} = n \frac{\lambda_n^r}{r!} e^{-\lambda_n}\left(1 + O\left(\frac{(r+1)^2}{n}\right)\right).$$

Hence, (8.5) is equivalent to

(8.7) $$n^{-1/2}(u_r(n) - \pi_r(\lambda_n)n) \xrightarrow{d} U_r \sim N(0, \varphi_{rr}).$$

The proof extends immediately to finite linear combinations of $u_r(n)$, which shows joint convergence in (8.5) and (8.7) for all $r \geq 0$; the covariances of the limits are given by

$$(8.8) \quad \varphi_{rs} := \operatorname{Cov}(U_r, U_s) = \pi_r(\lambda)\pi_s(\lambda)\bigg(\frac{(r-\lambda)(s-\lambda)}{\lambda} - 1\bigg) + \pi_r(\lambda)\delta_{rs}.$$

Each of $b_n$, $h_n$ and $l_n$ can be written as $\sum_{r=0}^{\infty} a_r(p)u_r(n)/n$, where $a_r(p)$ are some continuous functions on $[0,1]$, not depending on $n$ and such that $a_r(p) = O(r)$. We claim that, for any such $a_r(p)$, we have in $C[0,1]$ (with the usual uniform topology)

$$(8.9) \quad n^{-1/2}\sum_{r=0}^{\infty} a_r(p)(u_r(n) - \mathbb{E}\, u_r(n)) \xrightarrow{\mathrm{d}} \sum_{r=0}^{\infty} a_r(p)U_r.$$

Indeed, by the joint convergence in (8.5), this holds for the partial sums $\sum_{r=0}^{R}$ for any finite $R$. Moreover, it follows easily from the exact formula for $\operatorname{Var} u_r(n)$ (see [18], Example 6.35) that, for any given $M$ such that $\sup_n \lambda_n \leq M$,

$$(8.10) \quad n^{-1}\operatorname{Var} u_r(n) = O(M^r/r!)$$

uniformly in $r \geq 0$, and it is then routine to let $R \to \infty$ to obtain (8.9); see [2], Theorem 4.2. Further, using (8.6) again, it follows that, in $C[0,1]$,

$$(8.11) \quad n^{-1/2}\sum_{r=0}^{\infty} a_r(p)(u_r(n) - \pi_r(\lambda_n)n) \xrightarrow{\mathrm{d}} \sum_{r=0}^{\infty} a_r(p)U_r.$$

It is shown in [15] that the above results for $G(n, \lambda_n/n)$ combined with a simple monotonicity argument show that (8.7) holds, jointly for all $r \geq 0$, for $G(n,m)$ too, except that in this case $U_r$, $r \geq 0$, have a different joint Gaussian distribution with covariances

$$(8.12) \quad \varphi_{rs} := \operatorname{Cov}(U_r, U_s) = \pi_r(\lambda)\pi_s(\lambda)\bigg(-\frac{(r-\lambda)(s-\lambda)}{\lambda} - 1\bigg) + \pi_r(\lambda)\delta_{rs}.$$

(We still have $\mathbb{E}\, U_r = 0$.) Furthermore, elementary calculations, which we omit, show that the estimates (8.6) (ignoring the middle part) and (8.10) hold for $G(n,m)$ too. Hence, again by first considering finite sums $\sum_0^R$ and letting $R \to \infty$, (8.11) holds in $C[0,1]$ for $G(n,m)$ too.

Thus, (8.11) holds for both $G(n,\lambda_n/n)$ and $G(n,m)$, although the two cases have different $U_r$, which implies that (8.1)–(8.3) hold jointly, with $U_B, U_H, U_L$ given as linear combinations of $U_r$, $1 \leq r < \infty$. More precisely, let

$$T_r(p) := \sum_{l=r}^{\infty} \beta_{lr}(p)U_l,$$

$$Q_j(p) := \sum_{r=j}^{\infty} T_r(p)$$

and, so as to preserve the analogy with (5.33)–(5.35), let
$$U_W(p) := p^2 \sum_{r=1}^{\infty} rU_r.$$

Then
$$U_B(p) := Q_k(p),$$
$$U_H(p) := kQ_k(p) + \sum_{r=k+1}^{\infty} Q_r(p),$$
$$U_L(p) := U_W(p) - U_H(p).$$

It is easy to see that $U_W = 0$ for $G(n,m)$, and thus, $U_L = -U_H$, since $l_n(p) + h_n(p)$ is deterministic by (2.7) when the number of edges is, or by (8.9) with $a_r = r$.

Define, with $\varphi_{rs}$ given by (8.8) for $G(n, \lambda_n/n)$, and by (8.12) for $G(n,m)$,

$$(8.13) \quad \psi_{ij}(p) := \mathrm{Cov}(T_i(p), T_j(p)) = \sum_{l=i}^{\infty} \sum_{r=j}^{\infty} \beta_{li}(p)\beta_{rj}(p)\varphi_{lr},$$

$$(8.14) \quad \psi_{iW}(p) := \mathrm{Cov}(T_i(p), U_W(p)) = p^2 \sum_{l=i}^{\infty} \sum_{r=1}^{\infty} \beta_{li}(p) r \varphi_{lr}$$

[$\psi_{iW}(p) = 0$ for $G(n,m)$ because then $U_W = 0$], and further,

$$(8.15) \quad \sigma_{ij}^*(p) := \mathrm{Cov}(Q_i(p), Q_j(p)) = \sum_{l=i}^{\infty} \sum_{r=j}^{\infty} \psi_{lr}(p),$$

$$(8.16) \quad \sigma_{iW}^*(p) := \mathrm{Cov}(Q_i(p), U_W(p)) = \sum_{l=i}^{\infty} \psi_{lW}(p),$$

$$(8.17) \quad \sigma_{WW}^*(p) := \mathrm{Var}(U_W(p)) = \begin{cases} 2p^4\lambda, & \text{for } G(n, \lambda_n/n), \\ 0, & \text{for } G(n,m). \end{cases}$$

Then the expressions for $\sigma^*$ are analogous to (5.42)–(5.47), which by (8.15) and (8.16) and changes of summation order simplifies to

$$(8.18) \quad \sigma_{BB}^*(x) = \sum_{i=k}^{\infty} \sum_{j=k}^{\infty} \psi_{ij}(p),$$

$$(8.19) \quad \sigma_{BH}^*(x) = \sum_{i=k}^{\infty} \sum_{j=k}^{\infty} i\psi_{ij}(p),$$

$$(8.20) \quad \sigma_{HH}^*(x) = \sum_{i=k}^{\infty} \sum_{j=k}^{\infty} ij\psi_{ij}(p),$$

$$\sigma^*_{BL}(x) = \sum_{i=k}^{\infty} \psi_{iW}(p) - \sigma^*_{BH}(x), \tag{8.21}$$

$$\sigma^*_{HL}(x) = \sum_{i=k}^{\infty} i\psi_{iW}(p) - \sigma^*_{HH}(x), \tag{8.22}$$

$$\sigma^*_{LL}(x) = \sigma^*_{WW}(x) - 2\sum_{i=k}^{\infty} i\psi_{iW}(p) + \sigma^*_{HH}(x). \tag{8.23}$$

Note that $\psi_{iW}(p)$, $\sigma^*_{iW}(p)$ and $\sigma^*_{WW}(p)$ vanish for $G(n,m)$. □

We will now use Theorems 8.3 and 8.5 to prove Theorems 1.2 and 1.3.

PROOF OF THEOREM 1.2. By Lemma 8.4 and Theorem 8.3, Theorem 3.4 holds, with $b_n, h_n, l_n$ and $\hat{p}_n$ random. Furthermore, Theorem 8.5 applies, and since $b_n, h_n, l_n$ are independent of $Z_B, Z_H, Z_L$ for every $n$ by Theorem 8.3, $(U_B, U_H, U_L)$ and $(Z_B, Z_H, Z_L)$ are independent. By Lemma 8.2, we may assume that Condition 2.1 holds a.s., and by including further variables in the application of the Skorokhod coupling theorem in the proof of Lemma 8.2, we may further assume that (8.1), (8.2), (8.3), $\hat{p}_n \to \hat{p}$ and (3.10) hold a.s. (This trick is not essential, but we find it convenient to argue pointwise in the probability space, i.e., for each realization of the family of random processes. Note that the $O$ and $o$ terms that appear in this proof and the following are not assumed to be uniform over all points in the probability space; the implicit constants may thus be random, but they do not depend on $n$.)

By (2.13), for $0 < p \leq 1$,

$$l(p) = 0 \iff p = \psi_{k-1}(\lambda p) \iff \lambda = \frac{\lambda p}{\psi_{k-1}(\lambda p)}. \tag{8.24}$$

Since $\lambda > c_k := \inf_{\mu>0} \mu/\psi_{k-1}(\mu)$, the equation $\lambda = \mu/\psi_{k-1}(\mu)$ has at least one solution $\mu > 0$; see [17, Lemma 7.1]. We have defined $\hat{\mu} := \mu_k(\lambda)$ to be the largest solution, and thus, $p = \hat{\mu}/\lambda$ is the largest solution of (8.24). Since $\hat{\mu}/\lambda = \psi_{k-1}(\hat{\mu}) < 1$, we have $\hat{p} = \hat{\mu}/\lambda$ with $0 < \hat{p} < 1$.

Since $\hat{p} = \psi_{k-1}(\lambda\hat{p}) = \psi_{k-1}(\hat{\mu})$, (2.13) yields

$$\alpha := l'(\hat{p}) = \lambda\hat{p}(1 - \lambda\psi'_{k-1}(\lambda\hat{p})) = \lambda(\psi_{k-1}(\hat{\mu}) - \hat{\mu}\pi_{k-2}(\hat{\mu})). \tag{8.25}$$

By [17], Lemma 7.2 and its proof, for $k \geq 3$, and a simple calculation for $k = 2$, $\alpha := l'(\hat{p}) > 0$.

Recall that, as in Theorems 3.4 and 8.3, $\hat{p}_n$ is the largest zero in $[0,1]$ of $l_n$, which is random, and that (8.3) holds a.s. We define, for $n$ so large that $\lambda_n > c_k$, the nonrandom $\hat{p}^*_n$ as the largest value in $[0,1]$ such that $l_{\text{Po}(\lambda_n)}(\hat{p}^*_n) = 0$.

Applying the above argument to $\lambda_n$, $\hat{p}_n^* = \hat{\mu}_n/\lambda_n$. Furthermore, by (8.3), a.s.,

(8.26) $\qquad l_n(\hat{p}_n^*) - l_n(\hat{p}_n) = l_n(\hat{p}_n^*) = n^{-1/2} U_L(\hat{p}_n^*) + o(n^{-1/2}).$

Since $l_{\text{Po}(\lambda_n)} \to l = l_{\text{Po}(\lambda)}$ uniformly together with its derivatives, by (2.13) (or using the same proof as for Lemma 2.3), it follows easily (cf. the proof of Theorem 3.4) that $\hat{p}_n^* \to \hat{p}$; moreover, a.s., $\hat{p}_n \to \hat{p}$ by Theorem 3.4, and $l_n' \to l'$ uniformly by Lemma 2.3. Hence, Taylor's formula yields a.s.

$$l_n(\hat{p}_n^*) - l_n(\hat{p}_n) = (\hat{p}_n^* - \hat{p}_n)(l'(\hat{p}) + o(1)),$$

which, together with (8.26), yields

$$\hat{p}_n - \hat{p}_n^* = -n^{-1/2}\alpha^{-1} U_L(\hat{p}_n^*) + o(n^{-1/2}).$$

A Taylor expansion of $b_n$ now yields a.s., using Lemma 2.3 and (8.1),

$$b_n(\hat{p}_n) = b_n(\hat{p}_n^*) + b'(\hat{p})(\hat{p}_n - \hat{p}_n^*) + o(n^{-1/2})$$
$$= b_{\text{Po}(\lambda_n)}(\hat{p}_n^*) + n^{-1/2} U_B(\hat{p}) - n^{-1/2}\alpha^{-1} b'(\hat{p}) U_L(\hat{p}) + o(n^{-1/2}),$$

and similarly for $h_n(\hat{p}_n)$. Combined with (3.10), this shows that, with $W_\nu = Z_\nu(\hat{t}) + U_\nu(\hat{p})$, $\nu \in \{B, H, L\}$,

$$n^{-1/2}(v(\text{Core}_k) - b_{\text{Po}(\lambda_n)}(\hat{p}_n^*)n, e(\text{Core}_k) - \tfrac{1}{2} h_{\text{Po}(\lambda_n)}(\hat{p}_n^*)n)$$
$$\overset{\text{d}}{\longrightarrow} (W_B - \alpha^{-1} b'(\hat{p}) W_L, \tfrac{1}{2} W_H - \tfrac{1}{2}\alpha^{-1} h'(\hat{p}) W_L).$$

Since $\lambda \hat{p} = \hat{\mu}$, Lemma 2.4 yields

(8.27) $\qquad b(\hat{p}) = \psi_k(\hat{\mu}) \quad \text{and} \quad h(\hat{p}) = \hat{\mu} \psi_{k-1}(\hat{\mu}).$

Similarly, for $\lambda_n$, $b_{\text{Po}(\lambda_n)}(\hat{p}_n^*) = \psi_k(\lambda_n \hat{p}_n^*) = \psi_k(\hat{\mu}_n)$ and $h_{\text{Po}(\lambda_n)}(\hat{p}_n^*) = \hat{\mu}_n \times \psi_{k-1}(\hat{\mu}_n)$. Moreover, by Lemma 2.4 again,

(8.28) $\quad b'(\hat{p}) = \lambda \psi_k'(\lambda \hat{p}) = \lambda \pi_{k-1}(\hat{\mu}),$

(8.29) $\quad h'(\hat{p}) = \lambda(\psi_{k-1}(\hat{\mu}) + \hat{\mu} \psi_{k-1}'(\hat{\mu})) = \lambda(\psi_{k-1}(\hat{\mu}) + \hat{\mu} \pi_{k-2}(\hat{\mu})).$

Thus, by (8.25), $\alpha^{-1} b'(\hat{p}) = a_v$ and $\alpha^{-1} h'(\hat{p}) = a_e$ given in (1.3) and (1.4).

Hence, the result follows with $Z_v := W_B(\hat{t}) - a_v W_L(\hat{t})$ and $Z_e := \tfrac{1}{2} W_H(\hat{t}) - \tfrac{1}{2} a_e W_L(\hat{t})$. Since $\text{Cov}(W_\nu, W_\varkappa) = \text{Cov}(Z_\nu(\hat{t}), Z_\varkappa(\hat{t})) + \text{Cov}(U_\nu(\hat{p}), U_\varkappa(\hat{p})) = \hat{\hat{\sigma}}_{\nu\varkappa}$ given in (1.5), it follows that (1.6)–(1.8) hold. Since the distribution of $(Z_B(\hat{t}), Z_H(\hat{t}), Z_L(\hat{t}))$ is nonsingular by Theorem 3.1, so is the distribution of $(W_B, W_H, W_L)$, and hence also the distribution of $(Z_v, Z_e)$. $\square$

PROOF OF THEOREM 1.3. This is similar to the proof of Theorem 1.2. By Lemma 8.4 and Theorem 8.3, Theorem 3.5 holds, with the modifications given in Theorem 8.3, and we have $D \sim \text{Po}(\lambda)$ with $\lambda = c_k$.

We begin by checking the conditions on the function $l = l_{\text{Po}(\lambda)}$ in Theorem 3.5. By the definition of $\lambda = c_k$, $\mu/\psi_{k-1}(\mu) \geq \lambda$ for all $\mu > 0$, and equality holds only for $\mu = \hat{\mu} := \mu_k(c_k)$; see [17], Lemma 7.2. It follows by (2.13) that $l(p) \geq 0$ for $0 < p \leq 1$ with equality only for $p = \hat{p} := \hat{\mu}/\lambda = \psi_{k-1}(\hat{\mu})$. Since $\hat{\mu} > 0$, we have $0 < \hat{p} < 1$.

Note that $\hat{\mu}$ is the maximum point of $\psi_{k-1}(\mu)/\mu$, and thus, $\frac{d}{d\mu}(\psi_{k-1}(\mu)/\mu) = 0$ for $\mu = \hat{\mu}$, which yields

$$(8.30) \qquad \psi'_{k-1}(\hat{\mu}) = \frac{\psi_{k-1}(\hat{\mu})}{\hat{\mu}} = \frac{\hat{p}}{\hat{\mu}} = \frac{1}{c_k}.$$

Since $\hat{p}$ is a minimum point of $l$, $l'(\hat{p}) = 0$. Further, differentiation of (2.13) [using (8.30)] yields

$$\beta := l''(\hat{p}) = \lambda \hat{p}(-\lambda^2 \psi''_{k-1}(\lambda \hat{p})) = -\lambda^2 \hat{\mu} \psi''_{k-1}(\hat{\mu})$$
$$= \lambda^2 \hat{\mu}(\pi_{k-2}(\hat{\mu}) - \pi_{k-3}(\hat{\mu})) = \lambda^2 (\hat{\mu} - k + 2)\pi_{k-2}(\hat{\mu}) = \lambda^2 \widehat{\beta},$$

with $\widehat{\beta}$ defined in (1.9). It is easily checked that, for $\mu = k - 2$,

$$\frac{\psi_{k-1}(\mu)}{\mu \psi'_{k-1}(\mu)} = \frac{\sum_{i=k-1}^{\infty} \pi_i(\mu)}{\mu \pi_{k-2}(\mu)} = \sum_{j=0}^{\infty} \frac{\mu^j}{(k-1)\cdots(k-1+j)}$$
$$< \sum_{j=0}^{\infty} \frac{1}{k-1}\left(\frac{k-2}{k}\right)^j = \frac{k}{2(k-1)} < 1, \qquad \square$$

and thus $\hat{\mu} \neq k - 2$ by (8.30). Hence, $\beta \neq 0$, and thus $\beta > 0$.

We may thus apply Theorem 3.5. As in the proof of Theorem 1.2, we may extend Lemma 8.2 and assume that all limits we know to hold in probability or distribution by the results above actually hold a.s.

Fix some $\delta_0 \in (0, \hat{p})$. As in Theorem 3.5, let $\bar{p}_n$ be the (now random) minimum point of $l_n$ in $[\delta_0, 1]$; similarly, let $\bar{p}_n^*$ be the (deterministic) minimum point of $l_{\text{Po}(\lambda_n)}$ in $[\delta_0, 1]$. By Theorem 3.5 and the assumption in the previous paragraph, $\bar{p}_n \to \hat{p}$ a.s. Further, an argument analogous to that used in the proof of Theorem 3.5 shows that $\bar{p}_n^* \to \hat{p}$.

Using (8.3), and the fact that $\bar{p}_n$ and $\bar{p}_n^*$ are minimum points,

$$(8.31) \qquad l_n(\bar{p}_n) \leq l_n(\bar{p}_n^*) = l_{\text{Po}(\lambda_n)}(\bar{p}_n^*) + n^{-1/2} U_L(\bar{p}_n^*) + o(n^{-1/2}),$$

$$(8.32) \quad l_{\text{Po}(\lambda_n)}(\bar{p}_n^*) \leq l_{\text{Po}(\lambda_n)}(\bar{p}_n) = l_n(\bar{p}_n) - n^{-1/2} U_L(\bar{p}_n) + o(n^{-1/2}).$$

Since $U_L$ is continuous, it follows that

$$(8.33) \qquad l_n(\bar{p}_n) = l_{\text{Po}(\lambda_n)}(\bar{p}_n^*) + n^{-1/2} U_L(\hat{p}) + o(n^{-1/2}).$$

For $x \geq 0$, we write $l(p, x) := l_{\text{Po}(x)}(p)$ [as given in (2.13), with $\lambda$ replaced by $x$]. Thus, $l(p) = l(p, \lambda)$ and $l_{\text{Po}(\lambda_n)}(p) = l(p, \lambda_n)$. By (2.13), using $l(\hat{p}, \lambda) =$

0 and (8.30),
$$\frac{\partial}{\partial x}l(\hat{p},\lambda) = \lambda\hat{p}(-\hat{p}\psi'_{k-1}(\lambda\hat{p})) = -\hat{\mu}\hat{p}\psi'_{k-1}(\hat{\mu}) = -\hat{p}^2.$$

Hence, using the simple Lemma 8.6 below,
$$l_{\text{Po}(\lambda_n)}(\bar{p}_n^*) = l(\bar{p}_n^*,\lambda_n) \sim -\hat{p}^2(\lambda_n - c_k).$$

Combining this with (8.33), we find that

(8.34) $\quad l_n(\bar{p}_n) = -\hat{p}^2(\lambda_n - c_k)(1+o(1)) + n^{-1/2}U_L(\hat{p}) + o(n^{-1/2}).$

Thus, the asymptotic behavior of $n^{1/2}l_n(\bar{p}_n)$ is the same as that of $-n^{1/2}(\lambda_n - c_k)$, within $O(1)$. Accordingly, each of the three cases in Theorem 1.3 is matched by the corresponding case in Theorem 3.5. In particular, this proves part (i).

In case (ii), (8.34) yields
$$n^{1/2}l_n(\bar{p}_n) \to \zeta := -\hat{p}^2\gamma + U_L(\hat{p}).$$

Taking $\sigma = \sigma_{LL}^{1/2}$ as in Theorem 3.5, we may assume $\sigma Z = -Z_L(\hat{t})$ with $\hat{t} := -\ln\hat{p}$ [indeed, this is how $Z$ was defined in the proof of Theorem 3.5, see the line following (6.20)]. Thus,
$$\sigma Z - \zeta = -Z_L(\hat{t}) - U_L(\hat{p}) + \hat{p}^2\gamma = -W_L + \hat{p}^2\gamma,$$

where $W_L := Z_L(\hat{t}) + U_L(\hat{p})$, just as in the proof of Theorem 1.2.

Furthermore, from (8.31) and (8.32),
$$l_{\text{Po}(\lambda_n)}(\bar{p}_n) \le l_{\text{Po}(\lambda_n)}(\bar{p}_n^*) + o(n^{-1/2}).$$

On the other hand, a Taylor expansion around the minimum point $\bar{p}_n^*$ yields

(8.35)
$$\begin{aligned}l_{\text{Po}(\lambda_n)}(\bar{p}_n) &= l_{\text{Po}(\lambda_n)}(\bar{p}_n^*) + \tfrac{1}{2}(l''_{\text{Po}(\lambda_n)}(\bar{p}_n^*) + o(1))(\bar{p}_n - \bar{p}_n^*)^2 \\ &= l_{\text{Po}(\lambda_n)}(\bar{p}_n^*) + (\beta/2 + o(1))(\bar{p}_n - \bar{p}_n^*)^2,\end{aligned}$$

and thus
$$\bar{p}_n - \bar{p}_n^* = O(l_{\text{Po}(\lambda_n)}(\bar{p}_n) - l_{\text{Po}(\lambda_n)}(\bar{p}_n^*))^{1/2} = o(n^{-1/4}).$$

Moreover, again by Lemma 8.6, $\bar{p}_n^* - \hat{p} = O(|\lambda_n - \lambda|) = O(n^{-1/2})$, so we see that $\bar{p}_n - \hat{p} = o(n^{-1/4})$. Consequently, from (8.1) [since $b = b_{\text{Po}(\lambda)}$], we obtain

$$\begin{aligned}b_n(\bar{p}_n) &= b_{\text{Po}(\lambda_n)}(\bar{p}_n) + O(n^{-1/2}) = b_{\text{Po}(\lambda)}(\hat{p}) + O(|\lambda_n - \lambda| + |\bar{p}_n - \hat{p}| + n^{-1/2}) \\ &= b(\hat{p}) + o(n^{-1/4}),\end{aligned}$$

and similarly for $h$. It follows that in the conclusion of Theorem 3.5(ii), the random $b_n(\bar{p}_n)$ and $h_n(\bar{p}_n)$ may be replaced with $b(\hat{p})$ and $h(\hat{p})$. Upon combining (2.13) with $l'(\hat{p}) = 0$, (8.29) simplifies to

$$(8.36) \qquad h'(\hat{p}) = 2\lambda\hat{p} - l'(\hat{p}) = 2\lambda\hat{p},$$

and Theorem 1.3(ii) follows.

In case (iii), $n^{-1/2} = o(\lambda_n - c_k)$, so (8.34) yields

$$(8.37) \qquad -l_n(\bar{p}_n) = \hat{p}^2(\lambda_n - c_k)(1 + o(1)).$$

As in the proof of Theorem 1.2, let $\hat{p}_n^* = \hat{\mu}_n/\lambda_n$ be the largest zero in $(0, 1]$ of $l_{\mathrm{Po}(\lambda_n)}$. By (8.33) and (8.37),

$$l_{\mathrm{Po}(\lambda_n)}(\hat{p}_n^*) - l_{\mathrm{Po}(\lambda_n)}(\bar{p}_n^*) = -l_{\mathrm{Po}(\lambda_n)}(\bar{p}_n^*) = -l_n(\bar{p}_n) + O(n^{1/2}) \sim \hat{p}^2(\lambda_n - c_k),$$

so, using a Taylor expansion similar to (8.35),

$$\hat{p}_n^* - \bar{p}_n^* \sim (2/\beta)^{1/2}\hat{p}(\lambda_n - c_k)^{1/2}.$$

Hence, by Taylor again, since $l'_{\mathrm{Po}(\lambda_n)}(\bar{p}_n^*) = 0$ and $l''_{\mathrm{Po}(\lambda_n)}(\bar{p}_n^*) \to l''(\hat{p}) = \beta$,

$$(8.38) \qquad l'_{\mathrm{Po}(\lambda_n)}(\hat{p}_n^*) \sim (2\beta)^{1/2}\hat{p}(\lambda_n - c_k)^{1/2}.$$

Let $\hat{p}_n$ be as in Theorem 3.4 (now with $l_n$ random). By (8.3),

$$l_{\mathrm{Po}(\lambda_n)}(\hat{p}_n) - l_{\mathrm{Po}(\lambda_n)}(\hat{p}_n^*) = l_{\mathrm{Po}(\lambda_n)}(\hat{p}_n) = l_{\mathrm{Po}(\lambda_n)}(\hat{p}_n) - l_n(\hat{p}_n)$$
$$= -U_L(\hat{p})n^{-1/2} + o(n^{-1/2}).$$

Hence, for large $n$, using a Taylor expansion at $\hat{p}_n^*$ [as $l''_{\mathrm{Po}(\lambda_n)}$ is uniformly bounded],

$$l_{\mathrm{Po}(\lambda_n)}(\hat{p}_n^* - n^{-1/4}) = l_{\mathrm{Po}(\lambda_n)}(\hat{p}_n^*) - n^{-1/4}l'_{\mathrm{Po}(\lambda_n)}(\hat{p}_n^*) + O(n^{-1/2})$$
$$= l_{\mathrm{Po}(\lambda_n)}(\hat{p}_n) + O(n^{-1/2}) - n^{-1/4}l'_{\mathrm{Po}(\lambda_n)}(\hat{p}_n^*).$$

Then from (8.38), for large $n$,

$$l_{\mathrm{Po}(\lambda_n)}(\hat{p}_n^* - n^{-1/4}) = l_{\mathrm{Po}(\lambda_n)}(\hat{p}_n) + O(n^{-1/2}) - n^{-1/4}(2\beta)^{1/2}\hat{p}(\lambda_n - c_k)^{1/2}$$
$$< l_{\mathrm{Po}(\lambda_n)}(\hat{p}_n),$$

since $n^{-1/2} = o(n^{-1/4}(\lambda_n - c_k)^{1/2})$. Similarly, we see that $l_{\mathrm{Po}(\lambda_n)}(\hat{p}_n) < l_{\mathrm{Po}(\lambda_n)} \times (\hat{p}_n^* + n^{-1/4})$. Thus, for large $n$, $\hat{p}_n^* - n^{-1/4} < \hat{p}_n < \hat{p}_n^* + n^{-1/4}$ and, by the mean value theorem,

$$\hat{p}_n - \hat{p}_n^* = \frac{l_{\mathrm{Po}(\lambda_n)}(\hat{p}_n) - l_{\mathrm{Po}(\lambda_n)}(\hat{p}_n^*)}{(2\beta)^{1/2}\hat{p}(\lambda_n - c_k)^{1/2}(1 + o(1))}$$
$$= -(2\beta)^{-1/2}\hat{p}^{-1}(\lambda_n - c_k)^{-1/2}n^{-1/2}(U_L(\hat{p}) + o(1)).$$

Hence, by (8.1) and the mean value theorem,

$$\begin{aligned} b_n(\hat{p}_n) &= b_{\text{Po}(\lambda_n)}(\hat{p}_n) + O(n^{-1/2}) \\ &= b_{\text{Po}(\lambda_n)}(\hat{p}_n^*) + (b'(\hat{p}) + o(1))(\hat{p}_n - \hat{p}_n^*) + O(n^{-1/2}) \\ &= b_{\text{Po}(\lambda_n)}(\hat{p}_n^*) - (2\beta)^{-1/2}\hat{p}^{-1}(\lambda_n - c_k)^{-1/2}n^{-1/2}(b'(\hat{p})U_L(\hat{p}) + o(1)), \end{aligned}$$

and similarly for $h$. Consequently, we can replace the random $b_n(\hat{p}_n)$ and $h_n(\hat{p}_n)$ in (3.12) by the deterministic $b_{\text{Po}(\lambda_n)}(\hat{p}_n^*)$ and $h_{\text{Po}(\lambda_n)}(\hat{p}_n^*)$, provided we also replace $Z' = -(2\beta)^{-1/2}Z_L(\hat{t})$ by $Z' - (2\beta)^{-1/2}U_L(\hat{p}) = -(2\beta)^{-1/2}W_L$. The proof is completed by (8.27), (8.28), (8.36) and (8.37), taking $Z' := -(2\beta)^{-1/2}\hat{p}^{-1}\lambda W_L = -(2\widehat{\beta})^{-1/2}\hat{p}^{-1}W_L$.

LEMMA 8.6.  *Let $l(p, x)$ be a twice continuously differentiable function in a neighborhood of a point $(\hat{p}, \hat{x})$, and assume that $\frac{\partial l}{\partial p}(\hat{p}, \hat{x}) = 0$ and $\frac{\partial^2 l}{\partial p^2}(\hat{p}, \hat{x}) \neq 0$. Then there exist $\delta, \varepsilon > 0$ such that, if $|x - \hat{x}| < \delta$, then there is a unique $p(x)$ with $|p(x) - \hat{p}| < \varepsilon$ such that $\frac{\partial l}{\partial p}(p(x), x) = 0$. Moreover, as $x \to \hat{x}$, $p(x) - \hat{p} = O(|x - \hat{x}|)$ and*

$$(8.39) \qquad l(p(x), x) = (x - \hat{x})\left(\frac{\partial l}{\partial x}(\hat{p}, \hat{x}) + o(1)\right).$$

PROOF.  The existence of $p(x)$ follows by the implicit function theorem applied to $\partial l/\partial p$, which also shows that $p(x)$ is a differentiable function. Hence, $p(x) - \hat{p} = p(x) - p(\hat{x}) = O(|x - \hat{x}|)$. Furthermore, $x \mapsto l(p(x), x)$ is differentiable and

$$\frac{d}{dx}l(p(x), x)|_{x=\hat{x}} = \frac{\partial l}{\partial p}(\hat{p}, \hat{x})p'(x) + \frac{\partial l}{\partial x}l(\hat{p}, \hat{x}) = \frac{\partial l}{\partial x}l(\hat{p}, \hat{x}),$$

and (8.39) follows. □

PROOF OF THEOREM 1.4.  Fix $x \in (-\infty, \infty)$ and let $m = m(n) := \lfloor c_k n/2 + xn^{1/2}\rfloor$ and $\lambda_n = 2m/n$. Note that $\lambda_n \to c_k$ and $n^{1/2}(\lambda_n - c_k) \to 2x$. Then

$$n^{-1/2}\left(M - \frac{c_k}{2}n\right) \leq x \quad \iff \quad M \leq m(n)$$

$$\iff \quad G(n, m) \text{ has a nonempty } k\text{-core}.$$

By Theorem 1.3(ii), the probability of this converges to $\Phi(2x\hat{p}^2/\sigma)$. □

Just as before [see (3.1)–(3.3)], we denote $\check{b}_n(t) := b_n(e^{-t})$, $\check{h}_n(t) := h_n(e^{-t})$, $\check{l}_n(t) := l_n(e^{-t})$. Note that these functions now are random functions of $t$. [They depend on both $t$ and $(d_i)_1^n$.]

REMARK 8.7. In the proofs of Theorems 1.2 and 1.3 the randomness in the limit comes from $U_\nu$ and $Z_\nu$, $\nu \in \{B, H, L\}$, and they always end up in the combination $W_\nu = Z_\nu(\hat{t}) + U_\nu(\hat{p})$. Moreover, the proofs use arguments similar to the proofs of Theorems 3.4 and 3.5, and there is some repetition. It would be more satisfactory, in our opinion, if one could make a direct proof of Theorems 1.2 and 1.3 that combines $Z$ and $U$ to $W$ at an early stage.

Indeed, if one could prove that Theorem 3.1 holds for $G(n, (d_i)_1^n)$ too (see Remark 7.4 and Appendix), then Theorem 8.3 would extend to show also that Theorem 3.1 would hold for $G_n$, with $\check{b}_n$, $\check{h}_n$, $\check{l}_n$ defined by (3.1)–(3.3) random. For $G(n, \lambda_n/n)$ and $G(n, m)$ (assuming as always $m = \lambda_n n/2$ and $\lambda_n \to \lambda > 0$), Theorem 8.5 would then imply that Theorem 3.1 holds with $\check{b}_n$, $\check{h}_n$, $\check{l}_n$ replaced by the deterministic $b_{\text{Po}(\lambda_n)}(e^{-t})$, $b_{\text{Po}(\lambda_n)}(e^{-t})$, $b_{\text{Po}(\lambda_n)}(e^{-t})$, and $Z_\nu(t)$ replaced by $Z_\nu(t) + U_\nu(e^{-t})$. Theorems 1.2 and 1.3 then would follow by the same proofs as Theorems 3.4 and 3.5.

Even if we cannot (yet?) make the first step in this argument rigorous, it gives a strong intuitive motivation. To strengthen this intuition further, observe that we can construct a random multigraph $G_n^*$ as follows: first let $(d_i)_1^n$ have the distribution of the sequence of vertex degrees in $G(n, \lambda_n/n)$ [or $G(n, m)$], and then, given $(d_i)_1^n$, let $G_n^* = G^*(n, (d_i)_1^n)$. [Note that, if we were to take $G(n, (d_i)_1^n)$ instead, we would get back $G(n, \lambda_n/n)$ or $G(n, m)$.] For this random multigraph, the proof just outlined works fine, because Theorem 3.1 holds for $G^*(n, (d_i)_1^n)$. On the other hand, the proofs of Theorems 1.2 and 1.3 above go through verbatim for $G_n^*$ too. This shows that both proofs (when applicable) have to give the same result, and is another intuitive motivation for Theorems 1.2 and 1.3, as well as for the appearance of the variables $W_\nu$.

## APPENDIX

We certainly conjecture that Theorem 3.1 holds for $G(n, (d_i)_1^n)$ too, although we have failed to prove it; see Remark 7.4. Since this extension would give an alternative (and somewhat simpler) proof of Theorems 1.2 and 1.3 (see Remark 8.7), we will discuss another idea of how to prove it here. In principle, it seems possible to use this method to give a proof of Theorem 3.1 for $G(n, (d_i)_1^n)$, but verifying the conditions has turned out to be much harder than we expected. It seems intuitively obvious that they hold, but our attempts to check them rigorously have ended up in very technical estimates that we so far have not had the energy to complete. Perhaps a reader can find a simple argument.

The idea is to use the following extension of Proposition 4.1, which yields joint convergence of the $M_n$ together with some integer-valued processes $N_n$ that converge to a Poisson process. In our conceived application, each time

Proposition 4.1 is used in Section 5, we would use Proposition A.1 with $N_n(t)$ being the number of loops plus the number of multiple edges (excluding loops) that have been created up to time $t$. We would also have to count loops and multiple edges created after stopping (i.e., in the core), for example, by stopping the processes in $M_n$ but letting $N_n$ continue, with half-edges being paired by a suitable process until none are left. It is well known that then $N_n(\infty)$, the number of loops and multiple edges in $G^*(n,(d_i)_1^n)$, converges to a Poisson distribution (see, e.g., [16]), which fits well with the conclusions below. [Moreover, by speeding up the process at the end, we can assume that it is completed before some finite time $T$, so that $N_n(\infty) = N_n(T)$.]

A compensator $A_n$ satisfying (ii) below is easily expressed as an integral $A_n(t) = \int_0^t a_n(s)\,ds$, where $a_n(t)$ is the rate which with loops or multiple edges are created. For example, if the last white ball colored red happens to be in a bin with $j$ remaining white balls, loops are created with rate $r$; the rate for creating multiple edges is more complicated and depends on the numbers of edges already created with that bin (= vertex) as an endpoint, and the remaining number of balls at the other endpoints. The rate $a_n(t)$ is thus random and fluctuates rapidly as a function of $t$, but it seems almost obvious that a law of large numbers holds and that the integral $A_n(t)$ converges to a nonrandom function as in (iii) below; however, as said above, we have failed to find a completely rigorous proof.

PROPOSITION A.1. *Suppose that, in addition to the assumptions in Proposition 4.1, the following holds, for some right-continuous filtration of $\sigma$-fields $\mathcal{F}_n = (\mathcal{F}_n(t))_{t \geq 0}$ with respect to which $M_n$ is a martingale:*

(i) $N_n(t)$, $t \geq 0$, *is an increasing adapted integer-valued stochastic process with $N_n(0) = 0$ and all jumps equal to $+1$. [That is, $\Delta N_n(t) \in \{0,1\}$ for all $t \geq 0$.]*

(ii) $A_n(t)$, $t \geq 0$, *is a continuous stochastic process with $A_n(0) = 0$ such that $N_n(t) - A_n(t)$ is a martingale. (Thus, $A_n$ is the compensator of $N_n$.)*

(iii) *For each fixed $t \geq 0$, as $n \to \infty$, $A_n(t) \xrightarrow{\mathrm{p}} A(t)$, where $A(t)$ is a nonrandom continuous real-valued function.*

*Then $(M_n, N_n) \xrightarrow{\mathrm{d}} (M, N)$ as $n \to \infty$, in $D([0, \infty); \mathbb{R}^{q+1})$, where $M$ is as in Proposition 4.1 and:*

- $N$ *is a Poisson process with intensity $dA(t)$, that is, an increasing integer-valued stochastic process with independent increments such that $N(0) = 0$ and $N(t) - N(s) \sim \mathrm{Po}(A(t) - A(s))$ when $0 \leq s \leq t$;*
- $M$ *and $N$ are independent.*

*In particular, for any $T < \infty$, the conditional distribution $\mathcal{L}(M_n \mid N_n(T) = 0)$ converges to the distribution of $M$.*

It is easily seen that, under additional (weak) conditions, we may take $T = \infty$ too in the last statement; we omit the details. Note that $M_n \xrightarrow{\mathrm{d}} M$ by Proposition 4.1 and $N_n \xrightarrow{\mathrm{d}} N$ by [11], Theorem VIII.3.36; hence, only the joint convergence is new.

PROOF OF PROPOSITION A.1. Note that $(M, N)$ is a stochastic process with independent increments, and that $M$ is a martingale.

We will rely heavily on results in Jacod and Shiryaev [11], Chapter VIII. We therefore use the notation in [11], Chapters II and VIII. (This includes writing $n$ as a superscript rather than subscript.) Let $X^n := (M^n, N^n)$ and $X := (M, N)$. We will sometimes use subscripts $M$, $N$ or $X$ to indicate the process under consideration; thus, the characteristics ([11], II.2.6) of $X^n$ are denoted by $(B_X^n, C_X^n, \nu_X^n)$, and so on.

We choose a continuous truncation function $h_M : \mathbb{R}^q \to \mathbb{R}^q$. [Thus, $h_M(x) = x$ when $|x|$ is small and $h_M(x) = 0$ when $|x|$ is large.] Further, let $\psi : \mathbb{R} \to \mathbb{R}$ be a continuous function with $\psi(x) = 1$ when $|x| \le 1$ and $\psi(x) = 0$ when $|x| \ge 2$, and define $h : \mathbb{R}^{q+1} \to \mathbb{R}^{q+1}$ by

$$h(x, y) := (\psi(y) h_M(x), \psi(|x|) \psi(2y) y).$$

Then, $h$ is a continuous truncation function in $\mathbb{R}^{q+1}$. Note that, since $\Delta N^n \in \{0, 1\}$, $h(\Delta X^n) = (h_M(\Delta M^n), 0)$. Consequently,

$$X^n(h)_t := X_t^n - \sum_{s \le t} (\Delta X_s^n - h(\Delta X_s^n)) = (M^n(h_M)_t, 0).$$

Similarly, since $M$ is continuous and $\Delta N \in \{0, 1\}$, $h(\Delta X) = 0$, and

$$X(h)_t := X_t - \sum_{s \le t} (\Delta X_s - h(\Delta X_s)) = (M_t, 0).$$

It follows immediately from the definitions [11], II.2.6 and II.2.16 that

$$B_X^n = (B_M^n, 0),$$
$$\widetilde{C}_X^{n,ij} = \begin{cases} \widetilde{C}_M^{n,ij}, & 1 \le i, j \le q, \\ 0, & i = q+1 \text{ or } j = q+1, \end{cases}$$

and similarly for $B_X$ and $\widetilde{C}_X$. Consequently, the conditions [Sup-$\beta_5$] and [$\gamma_5$-$\mathbb{R}_+$] (see [11], VIII.2.1–2.2) for $X^n$ and $X$ are equivalent to the corresponding conditions for $M^n$ and $M$.

The conditions [Sup-$\beta_5$] and [$\gamma_5$-$\mathbb{R}_+$], together with [$\hat{\delta}_5$-$\mathbb{R}_+$], hold for $M^n$ and $M$ by [11], Theorem VIII.3.12 and Lemma VIII.3.15 (note that [Var-$\beta_5$] implies [Sup-$\beta_5$]), since [11], VIII.3.14, follows from our (4.4) and [11], VIII.3.12(b)(ii), is (4.3); cf. the proof of [14], Proposition 9.1. Consequently, [Sup-$\beta_5$] and [$\gamma_5$-$\mathbb{R}_+$] hold for $X^n$ and $X$.

Next, let $g \in C_2(\mathbb{R}^{q+1})$ [11], VII.2.7; thus, $g$ is a bounded, continuous function $\mathbb{R}^{q+1} \to \mathbb{R}$ such that $g(x) = 0$ when $|x|$ is small. Let $g_0(y) := g(0, y)$; then $g_0 \in C_2(\mathbb{R})$.

Define, using the notation in [11], II.1.5, the process

$$\widehat{D}_t^n := (g * \nu_X^n)_t - (g_0 * \nu_N^n)_t. \tag{A.1}$$

Then $\widehat{D}_t^n$ is the predictable compensator of the process (with the definition [11], II.1.16)

$$D_t^n := (g * \mu^{X^n})_t - (g_0 * \mu^{N^n})_t := \sum_{s \leq t} g(\Delta X^n(s)) - \sum_{s \leq t} g_0(\Delta N^n(s))$$

$$= \sum_{s \leq t} (g(\Delta M^n(s), \Delta N^n(s)) - g(0, \Delta N^n(s))).$$

For every $\rho > 0$, there exists $\delta > 0$ such that, if $|x| \leq \delta$, then $g(x, 0) = 0$ and $|g(x, 1) - g(0, 1)| < \rho$. Let $K := \sup|g|$, $\Delta^* M^n(t) := \sup_{s \leq t} |\Delta M^n(s)|$ and $J_\delta^n(t) := \sum_{s \leq t} \mathbf{1}[|\Delta M^n(s)| > \delta]$. Then,

$$|D_t^n| \leq \sum_{s \leq t} (2K \mathbf{1}[|\Delta M^n(s)| > \delta] + \rho \Delta N^n(s)) = 2K J_\delta^n(t) + \rho N^n(t).$$

Thus, $\widehat{D}$ is $L$-dominated ([11], I.3.29) by $2K J_\delta^n(t) + \rho N^n(t)$, and thus also by $F_t := 2K J_\delta^n(t) + \rho A^n(t)$. Note that $F$ is increasing and that $\Delta F_t = 2K \times \mathbf{1}[|\Delta M^n(t)| > \delta]$. Consequently, the domination inequality ([11], Lemma I.3.30) yields, for all $\varepsilon, \eta > 0$,

$$\mathbb{P}(|\widehat{D}_t^n| \geq \varepsilon) \leq \frac{1}{\varepsilon}\left(\eta + \mathbb{E}\left(\sup_{s \leq t} \Delta F_s\right)\right) + \mathbb{P}(F_t \geq \eta)$$

$$\leq \frac{1}{\varepsilon}(\eta + 2K\,\mathbb{P}(J_\delta^n(t) \geq 1)) + \mathbb{P}(J_\delta^n(t) > 0) + \mathbb{P}(\rho A^n(t) \geq \eta)$$

$$= \frac{\eta}{\varepsilon} + \left(\frac{2K}{\varepsilon} + 1\right)\mathbb{P}(\Delta^* M^n(t) > \delta) + \mathbb{P}(\rho A^n(t) \geq \eta).$$

As remarked above, $[\hat{\delta}_5\text{-}\mathbb{R}_+]$ holds for $M^n$ and $M$, which by [11, VIII.3.5] yields $\Delta^* M^n(t) \xrightarrow{\text{P}} 0$ as $n \to \infty$. Consequently, for any $\varepsilon, \eta, \rho > 0$, recalling (iii),

$$\limsup_{n \to \infty} \mathbb{P}(|\widehat{D}_t^n| \geq \varepsilon) \leq \frac{\eta}{\varepsilon} + \mathbb{P}(\rho A(t) \geq \eta).$$

Letting first $\rho \to 0$ and then $\eta \to 0$, we see that $\mathbb{P}(|\widehat{D}_t^n| \geq \varepsilon) \to 0$, that is,

$$\widehat{D}_t^n \xrightarrow{\text{P}} 0. \tag{A.2}$$

Moreover, since $N^n$ is a point process, $\nu_N^n$ is the random measure $dA^n(t) \times \delta_1$, where $\delta_1$ denotes a point mass at 1. Similarly, $\Delta X_t = (0, \Delta N_t)$ equals $(0,0)$ or $(0,1)$, and thus, $\nu_X = dA(t) \times \delta_{(0,1)}$. Hence (cf. the proof of [11], Theorem VIII.3.36),

$$(g_0 * \nu_N^n)_t = g_0(1)A^n(t) \xrightarrow{\mathrm{p}} g_0(1)A(t) = g(0,1)A(t) = (g * \nu_X)_t.$$

Combining this with (A.1) and (A.2), we see that $(g * \nu_X^n)_t \xrightarrow{\mathrm{p}} (g * \nu_X)_t$, which verifies $[\delta_{5,2}\text{-}\mathbb{R}_+]$ and thus $[\delta_{5,1}\text{-}\mathbb{R}_+]$ for $X^n$ and $X$.

The result $X^n \xrightarrow{\mathrm{d}} X$ now follows by [11], Theorem VIII.2.17. □

**Acknowledgments.** This research was initiated during a visit by MJL to Uppsala University in April 2005, sponsored by the LSE Nordic Exchange Scheme.

| | |
|---|---|
| DEPARTMENT OF MATHEMATICS | DEPARTMENT OF MATHEMATICS |
| UPPSALA UNIVERSITY | LONDON SCHOOL OF ECONOMICS |
| PO BOX 480 | HOUGHTON STREET |
| SE-751 06 UPPSALA | LONDON WC2A 2AE |
| SWEDEN | UNITED KINGDOM |
| E-MAIL: svante.janson@math.uu.se | E-MAIL: m.j.luczak@lse.ac.uk |
| URL: http://www.math.uu.se/~svante/ | URL: http://www.lse.ac.uk/people/m.j.luczak@lse.ac.uk/ |